\documentclass{svmult}
\usepackage{amscd}
\usepackage{amssymb}
\usepackage{amsmath}
\usepackage{times}
 \usepackage[all] {xy}      

\newtheorem{theo}{Theorem}[section]
\newtheorem{prop}[theo]{Proposition}
\newtheorem{lem}[theo]{Lemma}
\newtheorem{cor}[theo]{Corollary}

\newtheorem{defi}[theo]{Definition}

\newtheorem{rema}[theo]{Remark}

\newcommand{\qeed}{\hfill\textrm{QED}\break\null}
\newenvironment{demo}{\noindent\textit{Proof.}~}{\qeed}

\def\bb{\mathbb}
\begin{document}
\title*{The special symplectic structure of binary cubics }
\author{Marcus J. Slupinski  \thanks{MJS was supported in part by the Math Research 
Institute, OSU}\and Robert J. Stanton\thanks{RJS was supported in part by NSF Grants 
DMS-0301133 and DMS-0701198} }
\institute{IRMA, Universit\'e Louis Pasteur (Strasbourg),
 7 rue Ren\'e Descartes,
F-67084 Strasbourg Cedex, France, \texttt{slupins@math.u-strasbg.fr}
\and Department of Mathematics, Ohio State 
University, 231 West 18th Avenue, 
Columbus OH 43210-1174, \texttt{stanton@math.ohio-state.edu}}
 
\maketitle
\centerline{\bf Abstract}

Let $k$ be a field of characteristic not $2$ or $3$.  Let $V$ be the $k$-space of binary cubic polynomials. The natural symplectic structure on $k^2$ promotes to a symplectic structure
$\omega$ on $V$ and from the natural symplectic action of $\textrm{Sl}(2,k)$  one obtains
the symplectic module $(V,\omega)$. We give a  complete analysis of this symplectic module from  the point of view of the associated moment map, its norm square $Q$ (essentially the classical discriminant) and the symplectic gradient of $Q$. Among the results are a symplectic derivation of the Cardano-Tartaglia formulas for the roots of a cubic, detailed parameters for all $\textrm{Sl}(2,k)$ and $\textrm{Gl}(2,k)$-orbits, in particular identifying a group structure on the set of $\textrm{Sl}(2,k)$-orbits of fixed nonzero discriminant, and a purely symplectic generalization of the classical Eisenstein syzygy for the covariants of a binary cubic. Such fine symplectic analysis is due to the special symplectic nature inherited from the ambient exceptional Lie algebra $\mathfrak G_2$. 
\section{Introduction}

Binary cubic polynomials have been studied since the nineteenth century, being the natural setting for a possible extension of the rich theory of binary quadratic forms. An historical summary of progress on this subject can be found in \cite{DiHst}, especially concerning results related to integral coefficients. While for a fixed binary cubic interesting questions remain open, e.g. its range in the integers, the number of solutions, etc., it is the structure of the space of all binary cubics that is the topic of this paper.

The space of binary cubics, we will take coefficients in a field, is an example of a prehomogeneous vector space under $\textrm{Gl}(2,k)$, and from this point of view has been thoroughly investigated. Beginning with the fundamental paper by Shintani \cite{Sh}, recast adelically in \cite{Wr}, an analysis of this pv sufficient to obtain the properties of the Sato-Shintani zeta function was done. Subsequently several descriptions of the orbit structure were obtained, in particular relating them to extensions of the coefficient field. A feature of this space, and some other prehomogeneous spaces, apparently never exploited is the existence of a symplectic structure which is preserved by the natural action of $\textrm{Sl}(2,k)$. 

The purpose of this paper is to expose the rich structure of the space of binary cubics when viewed as a symplectic module using the standard tools of equivariant symplectic geometry, viz. the moment map, its norm square, and its symplectic gradient i.e. the natural Hamiltonian vector field. The advantages are several: somewhat surprisingly, the techniques are universally applicable, with the only hypothesis that the fields not be of characteristic 2 or 3; there are explicit symplectic parameters for each orbit type (including the singular ones not studied previously) that are easily computed for any specific field; the computations are natural; we obtain new results for the space of binary cubics e.g. a group structure on orbits; we obtain ancient results for cubics, namely a symplectic derivation of the Cardano-Tartaglia formula for a root.

This paper arose as a test case to see the extent that we might push a more general project \cite{Sl-St} on Heisenberg graded Lie algebras. A symplectic module can be associated to every such graded Lie algebra  and in the case of the split Lie algebra $\mathfrak G_2$, this symplectic module  turns out to be isomorphic to the space of binary cubics with the $\textrm{Sl}(2,k)$ action mentioned above. Although our approach to binary cubics  is inspired by the general situation, in order to give an accessible and elementary presentation, we have made this paper essentially self-contained with only one or two results quoted without proof from \cite{Sl-St}.

The symplectic technology consists of the following. The moment map, $\mu$, maps the space of binary cubics ${\rm S}^3({k^2}^*)$ to the Lie algebra $\mathfrak{\mathfrak{sl}}(2,k)$ of $\textrm{Sl}(2,k)$. By means of the Killing form on  $\mathfrak{\mathfrak{sl}}(2,k)$ one obtains a scalar valued function $Q$ on ${\rm S}^3({k^2}^*)$, the norm square of $\mu$. Using the symplectic structure one constructs $\Psi$, the symplectic gradient of $Q$,  as the remaining piece of symplectic machinery. This symplectic module appears to be "special" in several ways, e.g. a consequence of our analysis is that all the  $\textrm{Sl}(2,k)$ orbits in ${\rm S}^3({k^2}^*)$ are co-isotropic (see \cite{Sl-St} for the general case).  Let us recall that  over the real numbers it has been shown that there is also a very strong link between special symplectic connections (see \cite{Ca-Schw}) and Heisenberg graded Lie algebras (called $2$-graded in \cite{Ca-Schw}).

Here is a more detailed overview of the paper. We will analyze each of the symplectic objects $\mu,Q, \Psi$ and determine for each of them their image, their fibre, the $\textrm{Sl}(2,k)$ orbits in each fibre, and explicit parameters and isotropy for each orbit type. 
This is all done with symplectic methods, so that furthermore we identify the symplectic geometric meaning of these fibres. For example, we show that the null space, $Z$, of $\mu$ is the set of multiples of cubes of  linear forms. As $\textrm{Sl}(2,k)$ preserves the null space, we obtain a decomposition into a collection of isomorphic Lagrangian orbits which we show are parametrized by $k^*/{k^*}^3$. Binary cubics whose moment lies in the nonzero nilpotent cone of $\mathfrak{\mathfrak{sl}}(2,k)$ turn out to be those  which contain a factor that is the square of a linear form. For these there is only one orbit, whose image under $\mu$  we characterize.  The pullback by means of $\mu$ of the natural symplectic structure on the image and the restriction of the symplectic form on ${\rm S}^3({k^2}^*)$ essentially coincide. The generic case is when the image of the moment map lies in the semisimple orbits of $\mathfrak{\mathfrak{sl}}(2,k)$. In this case the  $\textrm{Sl}(2,k)$ orbits are different from the  $\textrm{Gl}(2,k)$ orbits, in contrast to the earlier cases. Here each of the values of $Q$ in $k{^*}$ determine a collection of  $\textrm{Sl}(2,k)$ orbits for which we give symplectic parameters using a `sum of cubes' theorem. As a consequence we show that the orbits for a fixed nonzero value of $Q$ form a group (over $\mathbb{Z}$ see \cite{Bh}) which we explicitly identify. Interestingly, a binary cubic is in the orbit corresponding to the identity of this group if and only if it is reducible.  The set of binary cubics corresponding to a fixed nonzero value of $Q$ is not stable under $\textrm{Gl}(2,k)$. However the set of binary cubics for which the value of $Q$ belongs to a fixed nonzero  square class of $k$ is stable under $\textrm{Gl}(2,k)$ and we obtain an explicit parametrisation of all   $\textrm{Gl}(2,k)$ orbits on this set.

If the field of coefficients is specialized to say $\mathbb{C}$ then several of the results herein are known. For example, that the zero set of $Q$ is the tangent variety to $Z$, or that the generic orbit is the secant variety of $Z$ can be found in the complex algebraic geometric literature. For some other fields other results are in the literature. However, the use of symplectic methods is new to all these cases and gives a unifying approach that seems to make transparent many classic results. For example, a careful analysis of $\mu$ and $\Psi$ in the generic case leads to a proof of the Cardano-Tartaglia formula for a root of a cubic.  As another application we conclude the paper with a symplectic generalization of the classical Eisenstein  syzygy for the covariants (compare to \cite{MoDio},\cite{HM}) of a binary cubic.  This is interesting because there is an analogue of this form of the Eisenstein syzygy for the symplectic module associated to any Heisenberg graded Lie algebra (\cite{Sl-St}). Finally, we remark that the symplectic methodology used in this paper could be used to understand binary cubics over the integers or more general rings.

We are very pleased to acknowledge the support of our respective institutions that made possible extended visits. To the gracious faculty of the Universit\'e Louis Pasteur goes a sincere merci beaucoup from RJS. In addition, RJS wants to acknowledge the support of Max Planck Institut, Bonn, for an extended stay during which some of this research was done.

\vskip .1cm
\section{Binary cubics as a symplectic space}

Let $k$ be a field such that $char(k)\not=2,3$. The vector space ${k^2}^*$ has a 
symplectic structure
$$
\Omega(ax+by,a'x+b'y)=ab'-ba'.
$$
Functorially one obtains a symplectic structure on the set of binary cubics 
$$
{\rm S}^3({k^2}^*)=\{  ax^3+3bx^2y+3cxy^2+dy^3:\, a,b,c,d\in k \}
.$$
 Explicitly, if
$P=ax^3+3bx^2y+3cxy^2+dy^3$
and
$P'=a'x^3+3b'x^2y+3c'xy^2+d'y^3,$
\begin{equation}\label{bincubsymp}
\omega(P,P') = ad'-da'-3bc'+3cb'.
\end{equation}
In particular, we have
\begin{equation}\label{sympcube}
\omega(P,(ex+fy)^3) = P(f,-e).
\end{equation}
Hence for $ex+fy\not=0$, 
\begin{equation}\label{divides}
(ex+fy) \mid\textrm{P} \iff \omega(P,(ex+fy)^3)=0.
\end{equation}
This indicates that one can use the symplectic form $\omega$ to study purely algebraic 
properties
of the space of binary cubics. More generally, the interplay of symplectic methods and 
the algebra of binary cubics will be the primary theme of this paper.

The group
$$
\textrm{Sl}(2,k)=\{\begin{pmatrix}\alpha&\beta\\
\gamma&\delta\end{pmatrix}:\,\alpha\delta-\beta\gamma=1\}
$$
acts on ${k^2}^*$ via the transpose inverse:
\begin{equation}\label{groupaction}
\begin{pmatrix}\alpha&\beta\\
\gamma&\delta\end{pmatrix}
\cdot x=
\delta x-\beta y,\quad
\begin{pmatrix}\alpha&\beta\\
\gamma&\delta\end{pmatrix}
\cdot y=
-\gamma x+\alpha y,
\end{equation}
and this action identifies $\textrm{Sl}(2,k)$ with the group of transformations of 
${k^2}^*$ that preserve the
symplectic form $\Omega$, i.e. $\textrm{Sp}({k^2}^*,\Omega)$. It follows that the 
functorial action of $\textrm{Sl}(2,k)$ on ${\rm S}^3({k^2}^*)$ preserves
the symplectic form $\omega$.  There is no kernel of this action thus 
$\textrm{Sl}(2,k)\hookrightarrow Sp({\rm S}^3({k^2}^*),\omega)$. 

The Lie algebra $\mathfrak{\mathfrak{sl}}(2,k)$ acts on ${k^2}^*$ via the negative 
transpose:
\begin{equation}\label{algebraaction}
\begin{pmatrix}
\alpha&\beta\\
\gamma&-\alpha
\end{pmatrix}
\cdot\,
 x=-\alpha x-\beta y,
 \quad
 \begin{pmatrix}
\alpha&\beta\\
\gamma&-\alpha
\end{pmatrix}
\cdot\,
 y=-\gamma x+\alpha y,
\end{equation}
which in terms of differential operators acting on polynomial functions on $k^2$ 
corresponds to the action
\begin{equation}\label{algactionasdop}
\begin{pmatrix}
\alpha&\beta\\
\gamma&-\alpha
\end{pmatrix}
\cdot\,f=
\alpha(-x\partial_xf+y\partial_yf)-\beta y\partial_xf-\gamma x\partial_yf.
\end{equation}
In particular, this gives
 the following action of $\mathfrak{\mathfrak{sl}}(2,k)$ on cubics:
\begin{align}
\nonumber
x^3&\mapsto -3\alpha x^3-3\beta x^2y\\
\nonumber
x^2y&\mapsto -\gamma x^3-\alpha  x^2y-2\beta xy^2\\
\nonumber
xy^2&\mapsto -2\gamma x^2y+\alpha xy^2-\beta y^3\\
\nonumber
y^3&\mapsto -3\gamma xy^2+3\alpha y^3.
\end{align}
\subsection{Symplectic covariants}
Among the basic tools of equivariant symplectic geometry are the moment map ($\mu$), 
its norm square ($Q$) and the symplectic gradient of $Q$ ($\Psi$). The symplectic 
structure on ${\rm S}^3({k^2}^*)$ is not generic as it is consistent with one 
inherited from an ambient Heisenberg graded Lie algebra, hence the description 
"special". In \cite{Sl-St} in the setting of Heisenberg graded Lie algebras we derive 
the fundamental properties of the basic symplectic objects as well as give 
explanations for normalizing constants, and identify characteristic features of these 
special symplectic structures. For the purposes of this paper the explicit formulae 
will suffice.

\begin{defi}
(i) The moment map $\mu:{\rm S}^3({k^2}^*)\rightarrow \mathfrak{\mathfrak{sl}}(2,k)$ 
here is
\begin{equation}\label{explicitmoment}
\mu(ax^3+3bx^2y+3cxy^2+dy^3)=
\begin{pmatrix}
ad-bc&2(bd-c^2)\\
2(b^2-ac)&-(ad-bc)
\end{pmatrix}.
\end{equation}
\noindent(ii) The cubic covariant  $\Psi:{\rm S}^3({k^2}^*)\rightarrow {\rm 
S}^3({k^2}^*) $ is given by
\begin{align}\label{psiformula}
\nonumber
\Psi(P)=\mu(P)\cdot P&=
(-3a\alpha-3b\gamma)x^3
+(-3a\beta-3b\alpha-6c\gamma)x^2y\\
&+(-6b\beta+3c\alpha-3d\gamma)xy^2
+(-3c\beta+3d\alpha)y^3
\end{align}
where $P=ax^3+3bx^2y+3cxy^2+dy^3$ and
$$
\begin{pmatrix}
\alpha&\beta\\
\gamma&-\alpha
\end{pmatrix}
=
\begin{pmatrix}
ad-bc&2(bd-c^2)\\
2(b^2-ac)&-(ad-bc)
\end{pmatrix}.
$$

\noindent(iii) The normalised quartic invariant  $Q_n:{\rm S}^3({k^2}^*)\rightarrow k$ is
\begin{equation}\label{normquart}
Q_n(P)=-\,{\rm det}\mu(P)= (a^2d^2-3b^2c^2-6abcd+4b^3d+4ac^3).
\end{equation}
\end{defi}

Notice that $Q_n(P)$ is a multiple of the classic discriminant of the polynomial $P$.

\begin{rema} The symmetric role of the coordinates  $x$ and $y$ is implemented by 
$$
\textrm{J}=
\begin{pmatrix}
0&-1\\
1&0
\end{pmatrix}
$$
which satisfies $\textrm{J}\cdot x= y$, $\textrm{J}\cdot y=-x$ and 
$$
\textrm{J}\cdot(ax^3+3bx^2y+3cxy^2+dy^3)=-dx^3+3cx^2y-3bxy^2+ay^3.
$$
From $\eqref{explicitmoment}$ it follows that $\mu(\textrm{J}\cdot P)$ is the cofactor 
matrix of $\mu(P)$.
\end{rema}
\begin{rema}The set of symplectic covariants $\omega,\mu,\Psi,Q,Q_n$ defined above is 
not the only choice possible for the purposes of this article. One could just as well 
use
$$
\omega_{\lambda}=\lambda\omega,\quad
\mu_{\lambda}=\lambda\mu,\quad\Psi_{\lambda}=\lambda\Psi,\quad
 Q_{\lambda}=\lambda^2Q$$
where $\lambda\in k^*$.
\end{rema}

The moment map  is characterised by the identity
\begin{equation}\label{charofmom}
Tr(\mu(P) \xi)=-\frac{1}{3}\omega(\xi\cdot P,P)\qquad\forall P\in{\rm 
S}^3({k^2}^*),\,\forall \xi\in \mathfrak{sl}(2,k),
\end{equation}
which specialized to $\xi = \mu(P)$ gives a characterization of $\Psi$
\begin{equation} Q(P) = 8\omega(P,\Psi(P)).
\end{equation}

From \eqref{charofmom} one gets that $\mu$ is $\textrm{Sl}(2,k)$-equivariant:
$$
\mu(g\cdot P)=g\mu(P) g^{-1}
\qquad\forall P\in{\rm S}^3({k^2}^*),\,\forall g\in \textrm{Sl}(2,k),
$$
and $\mathfrak{sl}(2,k)$-equivariant:
$$
d\mu_P(\xi\cdot P)=[\xi,\mu(P)]\qquad\forall P\in{\rm S}^3({k^2}^*),\,\forall \xi\in 
\mathfrak{sl}(2,k).
$$
Here, $d\mu_P(Q)=2B_{\mu}(P,Q)$ where $B_{\mu}:{\rm S}^3({k^2}^*)\times{\rm 
S}^3({k^2}^*)\rightarrow \mathfrak{sl}(2,k)$
is the unique symmetric bilinear map such that $\mu(P)=B_{\mu}(P,P)$.

From the $\textrm{Sl}(2,k)$ and $\mathfrak{sl}(2,k)$ equivariance of $\mu$ one obtains 
the $\textrm{Sl}(2,k)$ and $\mathfrak{sl}(2,k)$ equivariance of $\Psi,Q$ and $Q_n$.  
Several useful relations among $\mu,\Psi$ and $Q$ are derived in \cite{Sl-St}. The 
following involves a relation between vanishing sets of symplectic covariants.
\begin{prop}\label{vansets}
 Let $P$ be a binary cubic. Then
$$
\mu(P)=0\,\Rightarrow\,\Psi(P)=0\,\Rightarrow\, Q(P)=0.
$$
\end{prop}
\begin{demo} Since $\Psi(P)=\mu(P)\cdot P$, it is obvious that  
$\mu(P)=0\Rightarrow\Psi(P)=0$.
Suppose that $\Psi(P)=0.$ Then by equation \eqref{charofmom}
$$
Tr(\mu(P)^2)=-\frac{1}{3}\omega(\Psi(P),P)=0.
$$
 But $\mu(P)^2+{\rm det}\mu(P)Id=0$ by the Cayley-Hamilton theorem, so ${\rm 
det}\mu(P)=0$ and hence $Q(P)=0$.
\end{demo}

From the invariant theory point of view a covariant is an $\textrm{Sl}(2,k)$ invariant 
in $\rm{S}^*({\rm S}^3({k^2}^*))\otimes {\rm S}^\ast({k^2}^*)$. Concerning completeness of the 
symplectic invariants one has the classic result of Eisenstein \cite{Ei}.
\begin{prop}
(i) $\mu, \Psi $, $Q$ and the identity generate the $\textrm{Sl}(2,k)$ invariants in ${\rm 
S}^3({k^2}^*)\otimes {\rm S}^\ast({k^2}^*)$.
\medskip

\noindent(ii) The only relation among them viewed as functions on $k^2$ is 
$$ \Psi(P)(\cdot)^2 - 9Q_n(P) P(\cdot)^2 = -\frac{9}{2}\Omega_{k^2}(\mu(P)\cdot,\cdot)^3,$$
here $\Omega$ is extended by duality to $k^2 \times k^2$.
\end{prop}
\begin{demo} We shall give a symplectic proof of the relation (ii) in \S 3.
\end{demo}

\begin{rema} There are two interesting results related by a simple scaling to the Eisenstein result. Fix $P\in {\rm S}^3({k^2}^*)$ with $Q_n(P) \ne 0$. One can associate 
to $P$ a type of Clifford algebra, $Cliff_P$, and in \cite{Ha} it is shown that the 
center of $Cliff_P$ is the coordinate algebra of the genus one curve $X^2 - 27 Q_n(P) 
= Z^3.$ The other result arises from the observation that we could work over, say, 
$\mathbb{Z}$ instead of $k$. Then in \cite{Mo} Mordell showed that all integral 
solutions $(X,Y,Z)$ to $X^2 + kY^2 = Z^3$ with $(X,Z)=1$ are obtained from some 
$P\in{\rm S}^3({\mathbb{Q}^2}^*)$ with $Q_n(P) = -4k$ and evaluating $(ii)$ at a 
lattice point in $\mathbb Q^2$. We will not use these results in this paper but we 
will give a symplectic proof at another time. 
\end{rema}

\begin{rema} The Proposition gives a complete description of binary cubics from the 
point of view of $\textrm{Sl}(2,k)$ invariant theory. From the symplectic theory point 
of view, in \cite{Sl-St} we give characterizations of $\textrm{Sl}(2,k)$ as the 
subgroup of $Sp({\rm S}^3({k^2}^*),\omega)$ that preserves $Q(\cdot)$ and as the subgroup 
of $Sp({\rm S}^3({k^2}^*),\omega)$ that commutes with $\Psi$. 
\end{rema}

\subsection{The image of the moment map}
As $\mu:{\rm S}^3({k^2}^*)\rightarrow \mathfrak{\mathfrak{sl}}(2,k)$ is equivariant, 
the image of $\mu$ is a union of $\textrm{Sl}(2,k)$ invariant sets. Of course, the 
invariant functions on $\mathfrak{sl}(2,k)$ are generated by $det$. The following 
description of the orbits of $\textrm{Sl}(2,k)$ acting on level sets of $det$ uses the 
symplectic structure on ${k^2}^*$. Lacking any reference for this probably known 
result we include a proof. Subsequently, Paul Ponomarev brought to our attention the material in \cite{Bour} p.158-159 from which an alternate albeit non-symplectic proof can be extracted. 
\begin{prop}\label{sl2adjorb}
Let $\Delta\in k$ and set 
\begin{align}
\nonumber
\mathfrak{sl}(2,k)_{\Delta}&=\{X\in \mathfrak{sl}(2,k)\setminus\{0\}:\;{\rm 
det}X=\Delta\},\\
\nonumber
k_{\Delta}^*&=\{x\in k^*:\,\exists a,b\in k\text{ such that } x=a^2+b^2\Delta\}.
\end{align}
Then the orbits of $\textrm{Sl}(2,k)$ acting on $\mathfrak{sl}(2,k)_{\Delta}$ are in 
bijection with $k^*/k_{\Delta}^*$ under the map 
$\nu_{\Delta}:{\mathfrak{sl}(2,k)_{\Delta}}\rightarrow k^*/k_{\Delta}^*$ defined by
\begin{equation}\label{}
\nu_{\Delta}(X)=[\Omega(v,X\cdot v)]
\end{equation}
where $v$ is any element in ${k^2}^*$ which is not an eigenvector of $X$.
\end{prop}
\begin{demo} We make some preliminary remarks before proving the result. First we 
observe that the definition of $\nu_{\Delta}(X)$  is independent of choice of $v$. Indeed, given $v$ which is not an eigenvector of $X$, then $\{v,X\cdot v\}$ is a basis of ${k^2}^*$.
Given $w$ any other vector which is not an eigenvector then $w = a v + b X\cdot v$, 
and using Cayley-Hamilton we obtain that $[\Omega(v,X\cdot v)] = [\Omega(w,X\cdot w)]$.

Next, note that if $X\in \mathfrak{sl}(2,k)$ there exists $g\in{\textrm{Sl}}(2,k)$ and 
$\beta, \gamma\in k$ such that
$$
gXg^{-1}=
\begin{pmatrix}
0&\beta\\
\gamma&0
\end{pmatrix}.
$$
So to prove the result, we need only consider matrices in $\mathfrak{sl}(2,k)_{\Delta}$ of the form
$
X=\begin{pmatrix}
0&\beta\\
\gamma&0
\end{pmatrix}
$
with either $\beta$ or $\gamma$  nonzero. Since
$$
\begin{pmatrix}
0&1\\
-1&0
\end{pmatrix}
\begin{pmatrix}
0&\beta\\
\gamma&0
\end{pmatrix}
\begin{pmatrix}
0&1\\
-1&0
\end{pmatrix}^{-1}
=
\begin{pmatrix}
0&-\gamma\\
-\beta&0
\end{pmatrix},
$$
we can further suppose that $\gamma\not=0$. Then $x$ is not an eigenvector of $X$
and $\nu_{{\rm det}X}(X)=[\Omega(x,X\cdot x)]=[\Omega(x,\gamma x)]=[\gamma]$.

Suppose $\begin{pmatrix}
0&\beta\\
\gamma&0
\end{pmatrix}$ and $\begin{pmatrix}
0&\beta'\\
\gamma'&0
\end{pmatrix}$ in $\mathfrak{sl}(2,k)_{\Delta}$ have the same value of $\nu_{\Delta}$, i.e., $\beta\gamma=-\Delta=\beta'\gamma'$ and $[\gamma] = [\gamma']$.

Then there
exist $p,q$ in $k$ such that 
$\gamma'=(p^2+q^2{\rm det}X)\gamma$. Take as {\it Ansatz}
$$
\begin{pmatrix}
a&b\\
c&d
\end{pmatrix}=
\begin{pmatrix}
p&-q\frac {\Delta}{\gamma'}\\
\gamma q &p \frac{\gamma}{\gamma'}
\end{pmatrix}.
$$

Then 
$$
\begin{array}{l}
\text{ det }\begin{pmatrix}
a&b\\
c&d
\end{pmatrix} = p^2 \frac{\gamma}{\gamma'} + q^2 \Delta \frac{\gamma}{\gamma'} \\
= \frac{\gamma}{\gamma'} (p^2  + q^2 \Delta )\\
 =1.
\end{array}
$$

A routine computation shows that
$$
\begin{pmatrix}
a&b\\
c&d
\end{pmatrix}
\begin{pmatrix}
0&\beta'\\
\gamma'&0
\end{pmatrix}
\begin{pmatrix}
d&-b\\
-c&a
\end{pmatrix}
=
\begin{pmatrix}
0&\beta\\
\gamma&0
\end{pmatrix},
$$
and so $\nu_\Delta$ separates orbits.

To show that given $ \alpha \ne 0,$ there is an $X$ with $\text { det }X = 
\Delta$ and $\nu_\Delta (X) = [\alpha]$, take 
$$ X = \begin{pmatrix} 
0&-\frac{\Delta}{\alpha}\\
\alpha&0
\end{pmatrix}.
$$
Then $\text{ det }X = \Delta$ and $\nu_\Delta (X) = [\alpha]$. Finally, 
$\textrm{Sl}(2,k)$ invariance of $\nu_\Delta$ follows from the definition of 
$\nu_\Delta$.
\end{demo}

\begin{rema}
We make some elementary observations concerning the $\textrm{Sl}(2,k)$ adjoint orbits. 
 If $-\Delta\in {k^*}^2$, then $k_{\Delta}^*=k^*$ and there is only one orbit. If 
$\Delta=0$ then
$k_{\Delta}^*={k^*}^2$ and there is one nilpotent orbit for every element of 
$k^*/{k^*}^2$.  If $-\Delta\not\in{k^*}^2$ is nonzero, then $k_{\Delta}^*$ is  the 
set of  values in $k^*$ taken by the norm function associated to the quadratic 
extension
$k(\sqrt{-\Delta})$ or, equivalently, by the anisotropic quadratic form $x^2+\Delta 
y^2$
on $k^2$. It is well-known that this is a proper subgroup of $k^*$, at least in characteristic 0, (with thanks to P. Ponomarev for a discussion on characteristic p) and so in characteristic zero there are at 
least two orbits.
\end{rema}
\begin{rema}
Since $k^*/k_{\Delta}^*$ is a group, the Proposition puts a natural group structure 
on the set of orbits of ${ \textrm{Sl}(2,k)}$ acting on trace free matrices of fixed 
determinant.
Alternatively, $\mathfrak{sl}(2,k)$ can be ${ \textrm{Sl}(2,k)}$-equivariantly identified with  $S^2({k^2}^*)$, the space of binary quadratic forms, by
$$
X\quad
\longleftrightarrow \quad q_X(v)=\Omega(v,X\cdot v).
$$
By transport of structure, the Proposition then puts a natural group structure 
on  the set of orbits of ${\textrm{Sl}(2,k)}$ acting on binary quadratic forms of fixed discriminant. One can check that this is Gauss composition. In Theorems \ref{fixeddiscnonunitary} and \ref{fixeddisc} we will put a natural group structure on orbits of binary cubics with fixed nonzero
discriminant.
\end{rema}

The image of the moment map can be characterized as follows.
\begin{theo} \label{Immom}
Let $X\in \mathfrak{sl}(2,k)\setminus\{0\}$. Then
$$
X\in {\rm Im}\,\mu
\quad
\iff
\quad
\nu_{{\rm det}X}(X)=[2].
$$
\end{theo}
\begin{demo}
As before, we can suppose without loss of generality that
$
X=\begin{pmatrix}
0&\beta\\
\gamma&0
\end{pmatrix}
$
with say $\beta$  nonzero.

\noindent$(\Rightarrow):$ If $X=\mu(P)$ and $P=ax^3+3bx^2y+3cxy^2+dy^3$, we have
$$
\begin{array}{l}
ad-bc=0\\
2(bd-c^2)=\beta\\
2(b^2-ac)=\gamma.
\end{array}
$$
Hence $b\beta=d\gamma$ and  
$$
-\beta=2(c^2-d^2\frac{\gamma}{ \beta} )=2(c^2+(\frac{d}{ 
\beta})^2(-\beta\gamma))=2(c^2+(\frac{d}{\beta})^2{\rm det}X)
$$
so that $\nu_{{\rm det}X}(X)=[-\beta]=[2]$. 

\noindent $(\Leftarrow)$: Since  $\nu_{{\rm det}X}(X)=[-\beta]$ and by hypothesis 
$\nu_{{\rm det}X}(X)=[2]$, there
exist $p,q$ in $k$ such that 
$$-\beta=2(p^2+q^2{\rm det}X)=2(p^2-q^2\beta\gamma).
$$ 
If we set
$$
c=p,\quad a=\frac{\gamma}{\beta}p,\quad d=\beta q,\quad b=\gamma q
$$
and
$$
P=ax^3+3bx^2y+3cxy^2+dy^3,
$$
it is easily checked that 
$$
\mu(P)=
\begin{pmatrix}
ad-bc&2(bd-c^2)\\
2(b^2-ac)&-(ad-bc)
\end{pmatrix}=
\begin{pmatrix}
0&\beta\\
\gamma&0
\end{pmatrix}=
X.
$$
\end{demo}
\begin{rema} This result is a weak form of the Eisenstein identity. Indeed, if one 
cubes both sides of $\nu_{{\rm det}X}(X)=[2]$ and uses Gauss composition, one obtains 
the Eisenstein identity evaluated at a particular vector.
\end{rema}

\begin{rema} Varying the symplectic structure to $\omega_{\lambda}, \lambda \in k^*$ 
one can sweep out the other orbits with a moment map.
\end{rema}
\begin{cor} Let $P,P'$ be nonzero binary cubics such that $Q_n(P)=Q_n(P')$ and such that 
$\mu(P)$ and $\mu(P')$ are nonzero. Then there exists $g\in\textrm{Sl}(2,k)$ such that $g\cdot \mu(P)=\mu(P')$.
\end{cor}
\begin{demo} Since $Q_n(P)=Q_n(P')$, we have $\text { det }\mu(P)=\text { det }\mu(P')$. By the previous theorem,
$$
\nu_{\text { det }\mu(P)}(\mu(P))=\nu_{\text { det }\mu(P')}(\mu(P'))
$$
and the result follows from Proposition \ref{sl2adjorb}.
\end{demo}

\subsection{The image and fibres of $\Psi$}

\begin{prop} $P\in{\rm S}^3({k^2}^*)$ with $Q_n(P)\ne 0$  is in the image of $\Psi$ if and only if $9Q_n(P)$ is a cube in $k^*$.
\end {prop}
\begin{demo}\noindent$(\Rightarrow):$ Suppose that $P = \Psi(B)$. The key to the argument is a result from \cite{Sl-St} that is special to Heisenberg graded Lie algebras, namely a formula for $\Psi^2$. From this result one obtains $\Psi^2(B) = -(9Q_n(B))^2B$. On the other hand we have $\Psi^2(B) = \Psi(P)$. Hence $B= -(9Q_n(B))^{-2} \Psi(P)$. Applying $\Psi$ again and using that $\Psi$ is cubic we obtain $P=\Psi(B)=-\eta^3(9Q_n(P))^{2}P $ where $\eta = -(9Q_n(B))^{-2}$. So $(-\eta)^3=(9Q_n(P))^{-2}.$  Now $(-\eta (9Q_n(B))^{2})^3 =1$ so $(9Q_n(B))^{6} = (-\eta)^{-3} = (9Q_n(P))^{2}$. Thus we obtain $9Q_n(P) = (\pm 9Q_n(B))^3$.

\noindent$(\Leftarrow):$ Suppose $9Q_n(P) = \lambda^3$. Set $B= -\frac{1}{\lambda^2}\Psi(P)$. Then as above, $\Psi(B) = P$.
\end{demo}

\begin{cor} For $P\in{\rm S}^3({k^2}^*)$ with $9Q_n(P)\in {k^*}^3$ the fiber $\Psi^{-1} (P)$ consists of one element.
\end{cor}
\begin{demo} From the previous proof, if  $P = \Psi(B)$ then $B= -(9Q_n(B))^{-2} \Psi(P)$.  
\end{demo}

\begin{rema} We will see later that  a nonzero $P\in{\rm S}^3({k^2}^*)$ with $Q_n(P)= 0$ is in the image of $\Psi$ if and only if  $\mu(P)=0$ and $I_T(P)=[6]$ (cf Proposition \ref{muofpsi}). The fibre of $\Psi$ is then given by Proposition \ref{psifibresonN}.
\end{rema}

\section{Orbits and fibres}
\subsection{Symplectic covariants and triple roots}

One has the natural `algebraic'  condition
\begin{defi} 
$T=\{P\in{\rm S}^3({k^2}^*):\,P\not=0\text{ and }P\text{ has a triple root} \}$,
\end{defi}
\noindent and the natural `symplectic' condition
\begin{defi} 
$Z_{\mu}=\{P\in{\rm S}^3({k^2}^*):\,P\not=0\text{ and }\mu(P)=0\}$.
\end{defi}
The next proposition shows that the symplectic quantity $\mu$ detects the purely 
algebraic property of whether or not a binary cubic has a triple root.
\begin{prop}\label{tripleroot}
$T=Z_{\mu}$.
\end{prop}
\begin{demo} Let $P=ax^3+3bx^2y+3cxy^2+dy^3$. Then  $P\in Z_{\mu}$ iff
$\mu(P)=0$ iff $ad=bc,bd=c^2$ and $b^2=ac$. 

If $bc=0$, then $cbd=c^3=0$ and $b^3=acb=0$. Hence
$b=c=0$ and either $a=0$ or $d=0$. In the first case $P=dy^3$ and in the second 
$P=ax^3$.

If $bc\not=0$, then $a=\frac{b^2}{ c}$ and $d=\frac{c^2}{b}$ which means 
$P=\frac{1}{bc}(bx+cy)^3.$
\end{demo}

In order to determine the $\textrm{Sl}(2,k)$ orbit structure in the level set 
$Z_{\mu}=\mu^{-1}(0)\backslash \{0\}$ we need to construct an invariant that separates 
the orbits. We begin with the observation that the factorisation of $P\in T$ is not 
unique.
\begin{lem} \label{veronesescalng}
Let  $\lambda,\mu\in k^*$ and  $\phi,\psi\in{k^2}^*$  be such that
$
\lambda\phi^3=\mu\psi^3.
$
Then $\frac{\lambda}{\mu}$ is a cube and $\phi$ and $\psi$ are proportional.
\end{lem}
\begin{demo}
Unique factorisation.
 \end{demo}
This means the following (algebraic) definition makes sense.
\begin{defi} \label{defIT}
Define $I_T:T\rightarrow k^*/{k^*}^3$
by
$$
I_T(P)=[\lambda]_{k^*/{k^*}^3}
$$
where $P=\lambda\phi^3, \lambda\in k^*$ and $\phi\in{k^2}^*$.
\end{defi}

One can formulate the definition using symplectic methods. Given a non-zero
$\phi\in{k^2}^*$ there is a $g\in \textrm{Sl}(2,k)$ with 
$\Omega(\phi,g\cdot\phi) = 1$. If $P = \lambda \phi^3$ then 
\begin{equation}
\omega(P,(g\cdot\phi)^3) = \lambda\omega(\phi^3,(g\cdot\phi)^3) = 
\lambda\Omega(\phi,g\cdot\phi)^3 = \lambda.
\end{equation}
Thus $I_T(P) = [\omega(P,(g\cdot\phi)^3)].$

\begin{prop}
(i) Let $P_1,P_2\in T$. Then  
\begin{equation}
\textrm{Sl}(2,k)\cdot P_1=\textrm{Sl}(2,k)\cdot P_2 \iff I_T(P_1)=I_T(P_2).
\end{equation}
\medskip
\noindent(ii) The map $I_T$ induces a bijection of the space of orbits  
\begin{equation}Z_{\mu}/\textrm{Sl}(2,k)  \longleftrightarrow  k^*/{k^*}^3.
\end{equation}
\medskip
\noindent(iii) Let $P\in T$ and let $G_P=\{g\in \textrm{Sl}(2,k):\,g\cdot P=P\}$ be 
the isotropy subgroup of $P$. Then
$$
G_P=\{g\in \textrm{Sl}(2,k):\,\exists \mu\in k^*\text{ s.t. }g\cdot \phi=\mu\phi\text{ 
and }\mu^3=1\}
$$
 where $P=\lambda\phi^3, \lambda\in k^*$ and $\phi\in{k^2}^*$.
\end{prop}
\begin{demo} (i):
Suppose that $P_1=\lambda\phi^3$ and that there exists $g\in \textrm{Sl}(2,k)$ such 
that $g\cdot P_1=P_2$. 
Then $P_2=g\cdot (\lambda\phi^3)=\lambda(g\cdot\phi)^3$ and 
$I_T(P_2)=[\lambda]=I_T(P_1)$.

Conversely, suppose $P_1=\lambda_1\phi_1^3$, $P_2=\lambda_2\phi_2^3$ and 
$I_T(P_1)=I_T(P_2)$. The action
of  $\textrm{Sl}(2,k)$ on nonzero vectors of  ${k^2}^*$  is transitive so
we can find $g\in \textrm{Sl}(2,k)$ such that $g\cdot \phi_1=\phi_2$ and hence such 
that
$$
g\cdot P_1=\lambda_1\phi_2^3.
$$
Since $I_T(P_1)=I_T(P_2)$, there exists $\rho\in k$ such that 
$\lambda_1=\rho^3\lambda_2$ and 
$$
g\cdot P_1=\lambda_2(\rho\phi_2)^3.
$$
Choosing $h\in \textrm{Sl}(2,k)$ such that $h\cdot(\rho\phi_2)=\phi_2$, we have 
$(hg)\cdot P_1=P_2.$

\noindent(ii):  By (i), the map $I_T$ induces an injection of the space of orbits  of 
$\textrm{Sl}(2,k)$ acting on $T$ into  $k^*/{k^*}^3$. This is in fact a surjection 
since if $\lambda\in k^*$, $I_T(\lambda x^3)=[\lambda]$.

\noindent(iii):  This follows from unique factorisation.
\end{demo}
\begin{rema} Extending $\phi$ to a basis of ${k^2}^*$ we have the isomorphism
$$
G_P\cong\{\begin{pmatrix}
\mu&a\\
0&\frac{1}{ \mu}
\end{pmatrix}:\,
\mu\in k^*,\mu^3=1\text{ and }a\in k
\}.
$$
Consequently all the $\textrm{Sl}(2,k)$ orbits in $Z_\mu$ are isomorphic. Hence 
$Z_\mu$ is a smooth variety, and in \cite{Sl-St} we show that it is Lagrangian.

As the center of $\textrm{Gl}(2,k)$ acts on $Z_\mu$ by "cubes" it preserves $I_T$, and 
thus the $\textrm{Sl}(2,k)$ orbits in $Z_\mu$ are the same as the $\textrm{Gl}(2,k)$ orbits. From the point of view of algebraic 
groups, the result by Demazure \cite{Dema} characterizes $\textrm{Sl}(2,k)$ as the subgroup of the 
automorphisms of ${\rm S}^3({k^2}^*)$ that preserves $Z_\mu$.

\end{rema}

\subsection{Symplectic covariants and double roots}

In a similar way next we consider the `algebraic'  condition
\begin{defi} 
$D=\{P\in{\rm S}^3({k^2}^*):\,P\not=0\text{ and }P\text{ has a double root} \}$,
\end{defi}
\noindent and the `symplectic' condition
\begin{defi} 
$N_{\mu}=\{P\in{\rm S}^3({k^2}^*):\,P\not=0\text{ and }\mu(P)\text{ is nonzero 
nilpotent}\}$.
\end{defi}
  
Again it turns out that the  symplectic quantity $\mu$ detects the purely algebraic 
property of whether or not a binary cubic has a
double root.
\begin{theo} \label{doubleroot}
$D=N_{\mu}$.
\end{theo}
\begin{demo} The inclusion $D\subseteq N_{\mu}$ follows from the
\begin{lem} 
Let $P\in D$ and write  $P=(ex+fy)^2(rx+sy)$ with $ex+fy$ and $rx+sy$ independent. 
Then
$$
\mu(P)=\frac{2}{ 9}(es-fr)^2
\begin{pmatrix}
-ef&-f^2\\
e^2&ef
\end{pmatrix}.
$$
In particular,  ${\rm Ker}\,\mu(P)$ is spanned by the double root $ex+fy$.
\end{lem}
\begin{demo}
Straightforward calculation.
\end{demo}

To prove the inclusion $N_{\mu}\subseteq D$, suppose $\mu(P)$ is a nonzero nilpotent. 
Then ${\rm Ker}\,\mu(P)$ is one-dimensional, spanned by, say,
$v\in {k^2}^*$. Since $\textrm{Sl}(2,k)$ acts transitively on nonzero vectors in  
${k^2}^*$, there exists
$g\in \textrm{Sl}(2,k)$ such that 
$
g\cdot v=x$.
Then $\mu(g\cdot P)=g\mu(P)g^{-1}$ is nonzero nilpotent with kernel spanned by $x$. 
Let
$g\cdot P=ax^3+3bx^2y+3cxy^2+dy^3$. Then by the formulae \eqref{algactionasdop} and 
\eqref{explicitmoment},  the condition $\mu(g\cdot P) \cdot x=0$ is equivalent to the 
system
\begin{align}
\nonumber
ad-bc&=0\\
\nonumber
bd-c^2&=0.
\end{align}
If $(c,d)\not=(0,0)$ this implies there exists $\lambda,\nu\in k$ such that 
$(a,b)=\lambda(c,d)$ and
$(b,c)=\nu(c,d)$. Hence $c=\nu d, b=\nu^2 d, a=\nu^3 d$ and $\mu(g\cdot P)=0$ which is 
a contradiction.
Thus $c=d=0$ and $g\cdot P=ax^3+3bx^2y=x^2(ax+3by)$. We have  $b\not=0$ (otherwise 
$\mu(P)=0$) so
$x$ and $ax+by$ form a basis of $k^2$. Applying $g^{-1}$ to $g\cdot P=x^2(ax+3by)$ 
completes the proof.
\end{demo}

Again, in order to obtain parameters for the orbit structure of $N_\mu$ we need 
standard representatives. The factorisation of $P\in N_{\mu}$ given by Theorem 
\ref{doubleroot} is not unique. However we can use the symplectic form $\Omega$ on 
${k^2}^*$ to get a canonical form for $P$.
\begin{lem}
Let $P\in N_{\mu}$. There exists a unique basis $\{\phi,\xi\}$ of ${k^2}^*$ such that 
$P=\phi^2\xi$ and $\Omega(\phi,\xi)=1$.
\end{lem}
\begin{demo}  If $P\in N_{\mu}$ then $P$ has a double root by Theorem 
\ref{doubleroot}.  Fix a factorisation $P=\phi_1^2\xi_1$.  By unique factorisation, 
any other factorisation is of the form $P= \phi^2\xi$
where
$$
\phi=\lambda \phi_1,\quad\xi=\frac{1}{\lambda^2}\xi_1
$$
for some $\lambda\in k^*$. Then $\Omega(\phi,\xi)=1$ iff 
$\lambda=\Omega(\phi_1,\xi_1)$ and this proves the claim.
\end{demo}
\begin{prop} The group  $\textrm{Sl}(2,k)$ acts simply transitively on $N_{\mu}$. Consequently, $\textrm{Gl}(2,k)$ has one orbit on $N_{\mu}$.
\end{prop}
\begin{demo}
Let $P,Q\in N_{\mu}$ and write $P=\phi^2\xi$ and $Q={\phi'}^2\xi'$ with  
$\Omega(\phi,\xi)= \Omega(\phi',\xi')=1$. The element
$g$ of $GL(2,k)$ defined by $g\cdot \phi=\phi'$ and $g\cdot \xi=\xi'$ is clearly in 
$\textrm{Sl}(2,k)$,  satisfies $g\cdot P=Q$
and is the unique element of $\textrm{Sl}(2,k)$ sending $P$ to $Q$.
\end{demo}
\begin{rema} In \cite{Sl-St} when char $k = 0$ we show that $N_{\mu}$ is the tangent variety to $Z_\mu$.
\end{rema}

\begin{rema} From  Proposition \ref{tripleroot}, Theorem \ref{doubleroot}  and 
\eqref{normquart} we see that $Q_n(P)=0$ iff $P$
has a multiple root, which is consistent with the classic discriminant interpretation. Also, the open subset of double roots is isomorphic to 
$\textrm{Sl}(2,k)$. Consequently the variety $Q_n(P)=0$ is not smooth, but has 
singular set which is a union over $k^*/{k^*}^3$  of isomorphic Lagrangian $\textrm{Sl}(2,k)$-orbits. 
\end{rema}

\vskip0,1cm
\centerline{{\bf Image and Fibres of  $\mu:N_{\mu}\rightarrow {\rm \mathfrak{sl}}(2,k)$ }}
\vskip0,1cm

The image  of the moment map on 
$N_{\mu}$ is given by Theorem \ref{Immom}:
\begin{cor}
$\mu(N_{\mu})=\{X\in {\rm \mathfrak{sl}}(2,k)\setminus\{0\}:\,{\rm det}X=0\text{ and 
}\nu_0(X)=[2]\}$.
\end{cor}

Now we give two descriptions of the fibres of $\mu:N_{\mu}\rightarrow {\rm 
\mathfrak{sl}}(2,k)$: the first
symplectic, the second algebraic. Note that the fibres of the moment map are 
symplectic objects so it is not a priori clear that they have a purely algebraic 
description.

\begin{prop} \label{momfibnil}
Let  $P\in N_{\mu}$ and let $\phi\in{k^2}^*$ be a square factor of $P$.

\medskip
\noindent (a) 
$\mu^{-1}(\mu(P))=\{P+a\Psi(P):a\in k\}\cup\{-P+b\Psi(P):b\in k\}.$

\medskip
\noindent(b) $\mu^{-1}(\mu(P))=\{P+a\phi^3:\,a\in k\}\cup\{-P+b\phi^3:b\in k\}.$

\medskip
\noindent(c) The affine lines in (a) and (b) are disjoint.
\end{prop}
\begin{demo} Since $\textrm{Sl}(2,k)$ acts transitively on $N_{\mu}$ we can assume 
without loss of generality that 
$P=3x^2y$. Then by \eqref{explicitmoment} and \eqref{psiformula},
$$
\mu(3x^2y)=\begin{pmatrix}
0&0\\
2&0
\end{pmatrix},\quad
\Psi(3x^2y)=-6x^3.
$$
We want to find all   $Q\in{\rm S}^3({k^2}^*)$ such that
\begin{equation}\label{invimmunil}
\mu(Q)=
\begin{pmatrix}
0&0\\
2&0
\end{pmatrix}.
\end{equation}
By Theorem \ref{doubleroot},  a solution of this equation is of the form 
$Q=(ex+fy)^2(rx+sy)$  with $es-fr\not=0$.  Substituting back in \eqref{invimmunil}   
we get
$$
\frac{2}{9}(es-fr)^2
\begin{pmatrix}
-ef&-f^2\\
e^2&ef
\end{pmatrix}
=
\begin{pmatrix}
0&0\\
2&0
\end{pmatrix}
$$
from which it follows  that the set of solutions of equation \eqref{invimmunil} is:
$$
\{x^2(e^2rx+3y):\, e\in k^*,r\in k \}\cup\{x^2(e^2rx-3y):\, e\in k^*,r\in k \}.
$$
Since $P=3x^2y$ and $\Psi(P)=-6x^3$, this proves  (a), (b) and (c).
\end{demo}

The fibre of $\mu$ at  $\mu(P)$ is also the orbit through $P$ of the isotropy group of 
$\mu(P)$.
\begin{cor}
Let  $P\in N_{\mu}$ and let $G_{\mu(P)}=\{g\in 
\textrm{Sl}(2,k):\,g\mu(P)g^{-1}=\mu(P)\}$. Then
$\mu^{-1}(\mu(P))=G_{\mu(P)}\cdot P$.
\end{cor}
\begin{demo}  Since $\mu(P)$ is nilpotent nonzero, a simple calculation shows that
$$
G_{\mu(P)}=\{Id+a\mu(P):\,a\in k\}\cup\{-Id+b\mu(P):\,b\in k\}
$$
and the result follows from Proposition \ref{momfibnil}.
\end{demo}

It appears that $N_{\mu}$ is a regular contact variety. If one endows the nilpotent 
variety $\mathcal{N}$ in $\mathfrak{sl}(2,k)$ with the KKS symplectic structure, then 
$\mu : N_{\mu}\to \mathcal{N}$ is a prequantization of the image of $\mu$. 

 \vskip0,1cm
 \centerline{{\bf Image and Fibres of  $\Psi:N_{\mu}\rightarrow Z_{\mu}$ }}
\vskip0,1cm

We begin with some properties of $\Psi$.

 \begin{prop}\label{muofpsi}
 Let $P=\phi^2\xi$ with $\phi,\xi\in{k^2}^*$. Then:
 
\medskip
\noindent(i) $\mu(\Psi(P))=0$ ;
 
\medskip

\noindent(ii) $\phi^3$ divides $\Psi(P)$;
 
\medskip

\noindent(iii) $\Psi(P)=0$ iff $\mu(P)=0$ ;
 
\medskip

 \noindent(iv) $\Psi(P)\not=0$ \,$\Rightarrow$\,  $I_{T}(\Psi(P))=[6]_{k^*/{k^*}^3}$.
\end{prop}
\begin{demo}
Set $\phi=ex+fy$ and $\xi=rx+sy$. Then calculation gives 
 \begin{align}\label{doublerootmuPsi}
\nonumber
\mu(P)&=\frac{2}{9}(es-fr)^2
\begin{pmatrix}
-ef&-f^2\\
e^2&ef
\end{pmatrix},
\\
\Psi(P)&=-\frac{2}{9}(es-fr)^3(ex+fy)^3
\end{align}
and all parts of the proposition follow immediately from these formulae.
\end{demo}

\begin{cor} The image of $\Psi$ on $N_{\mu}$ is $Z_{\mu}[6]$.
\end{cor}
\begin{demo}
According to Proposition \ref{muofpsi}(iv), if $P\in N_{\mu}$ then $\Psi(P)\in 
Z_{\mu}$ and $I_{T}(\Psi(P))=[6]_{k^*/{k^*}^3}$.
Since $\Psi$ is $\textrm{Sl}(2,k)$-equivariant and $\textrm{Sl}(2,k)$ acts 
transitively on both $N_{\mu}$ and $Z_{\mu}[6]$, it is clear that 
$\Psi$ maps $N_{\mu}$ onto $Z_{\mu}[6]$.
\end{demo}

To  describe the fibres we need a symplectic 
characterization of the double root of a $P\in Z_{\mu}$. Recall that $ex+fy\ne 0$ is a 
root of $P$ iff $\omega(P,(ex+fy)^3 )= 0$. Analogous to this result we have

\begin{prop} Let $P$ be a binary cubic and $(ex+fy)\in {k^2}^*$ be nonzero. 
\begin{equation} (ex+fy)^2  \mid\textrm{P} \iff B_{\mu}(P,(ex+fy)^3)=0.
\end{equation}
\end{prop}
\begin{demo}
We begin with two remarks. First,
since $\textrm{Sl}(2,k)$ acts transitively on nonzero elements of ${k^2}^*$ and since 
$B_{\mu}$ and $\Psi$ are
$\textrm{Sl}(2,k)$-equivariant, we can assume without loss of generality that 
$ex+fy=x$.  Second, the formula  for
$B_{\mu}$ obtained by polarising \eqref{explicitmoment} is
\begin{equation}\label{explicitBmu}
B_{\mu}(P,P')=
\begin{pmatrix}
\frac{1}{2}(ad'+da'-bc'-cb')&(bd'+db')-cc'\\
bb'-(ac'+ca')&-\frac{1}{2}(ad'+da'-bc'-cb')
\end{pmatrix}
\end{equation}
if
$P=ax^3+3bx^2y+3cxy^2+dy^3$
and
$P'=a'x^3+3b'x^2y+3c'xy^2+d'y^3$.

\noindent Let $P=ax^3+3bx^2y+3cxy^2+dy^3$. Then
$$
B_{\mu}(P,x^3)=
\begin{pmatrix}
\frac{1}{2}d&0\\
-c&-\frac{1}{2}d
\end{pmatrix}
$$
and hence $x^2$ divides $P$ iff $c=d=0$ iff $B_{\mu}(P,x^3)=0.$
\end{demo}

Now since $\Psi$ maps $D$ to $T$ we expect a criterion involving $\Psi$ for $ex+fy\ne 
0$ to be a double root of $P$.

\begin{prop}Let $P$ be a binary cubic and $(ex+fy)\in {k^2}^*$ be nonzero.

 \medskip
\noindent(i) If $(ex+fy)^2$ divides $P$ then $\Psi(P)$ is proportional to $(ex+fy)^3$.

\medskip
\noindent(ii) If $\Psi(P)$ is a nonzero multiple of $(ex+fy)^3$ then $(ex+fy)^2$ 
divides $P$.

\medskip
\noindent(iii) $\{P\in {\rm S}^3({k^2}^*):\,B_{\mu}(P,(ex+fy)^3)=0\}$ is a Lagrangian 
subspace of ${\rm S}^3({k^2}^*)$.
\end{prop}
\begin{demo}
\noindent $(i)$: If $x^2$ divides $P$ then taking $e=1$ and $f=0$ in the formulae 
\eqref{doublerootmuPsi} we get $\Psi(P)=-\frac{2}{9}d^3x^3$. 

\noindent $(ii)$: If there exists $\lambda\in k^*$ such that 
$(ex+fy)^3=\frac{1}{\lambda}\Psi(P)$, we have
$$
B_{\mu}(P,(ex+fy)^3)=\frac{1}{\lambda}B_{\mu}(P,\mu(P)\cdot P).
$$
But $B_{\mu}(P,\mu(P)\cdot P)+B_{\mu}(\mu(P)\cdot P,P)=[\mu(P),\mu(P)]=0$ since 
$B_{\mu}$ is $\mathfrak{sl}(2,k)$-equivariant.
Hence $B_{\mu}(P,\mu(P)\cdot P)=0$ and $B_{\mu}(P,(ex+fy)^3)=0$  which implies by the 
previous result  that $(ex+fy)^2$ divides $P$.

\noindent $(iii)$:  Let $L=\{P\in {\rm S}^3({k^2}^*):\,B_{\mu}(P,(ex+fy)^3)=0\}$. As 
we saw in the proof above, the binary cubic $ax^3+3bx^2y+3cxy^2+dy^3$ is in  $L$ iff 
$c=d=0$ and hence $L$ is of dimension two. It follows from \eqref{bincubsymp} that 
$\omega(P,P')=0$ if $P,P'\in L$ and hence $L$ is Lagrangian.
\end{demo}

We can now give two descriptions of the fibres of $\Psi:N_{\mu}\rightarrow 
Z_{\mu}[6]$, the first
symplectic, the second algebraic. Again, as the fibres of $\Psi$ are symplectic 
objects it is not a priori clear that they have a purely algebraic description.

\begin{prop}  \label{psifibresonN}
Let $P\in N_{\mu}$ and let $\phi\in{k^2}^*$ be a square factor of $P$.

\medskip
(i)  $\Psi^{-1}(\Psi(P))=\{aP+b\Psi(P):\,a\in k^*,b\in k\}$.

\medskip
(ii)  $\Psi^{-1}(\Psi(P))=\{Q\in N_{\mu}:\,\phi^2\text{ divides }Q\}$.
\end{prop}

\vskip0.1cm
 \centerline{\bf Explicit factorisation of  $P$ when $Q_n(P)=0$}
\vskip0.1cm

From what has been done thus far we obtain readily

\begin{prop} Let $P=ax^3+3bx^2y+3cxy^2+dy^3$ be a nonzero binary cubic over a field 
$k$ such that ${\rm char}(k)\not=2,3$.

\medskip
\noindent(i) If $\mu(P)=0$ then $Q_n(P)=0$ and
$$
P=\left\{
\begin{array}{ll}
ax^3\text{ or } dy^3&\text{ if }bc=0,\\
\frac{1}{bc}(bx+cy)^3&\text{ if }bc\not=0.
\end{array}
\right.
$$
(ii) If $\mu(P)\not=0$ and $Q_n(P)=0$ then
$$
P=\left\{
\begin{array}{ll}
x^2(ax+3by)\text{ or }(3cx+d) y^2&\text{ if }ad-bc=0,\\
\bigl(-(b^2-ac)x+\frac{1}{2}(ad-bc)y\bigr)^2(\frac{a}{(b^2-ac)^2}x+\frac{4d}{(ad-bc)^2
}y)&\text{ if }ad-bc\not=0.
\end{array}
\right.
$$
\end{prop}

\subsection{ Symplectic covariants and  sums of coprime cubes}
We have seen that a $P$ with multiple roots corresponds to $Q_n(P) =0$. So we begin 
the study of $P$ with $Q_n(P)\ne 0$, in which case the  $\textrm{Sl}(2,k)$ orbits are not the same as the  $\textrm{Gl}(2,k)$ orbits. The values of the symplectic invariant $Q_n(P)$ 
will have much to say about the roots of $P$. We begin with the  `natural' condition
\begin{defi} 
${\cal O}_{[1]}=\{P\in{\rm S}^3({k^2}^*):\,Q_n(P) \text{ is a square in }{k}^*\}$.
\end{defi}

The relevant `algebraic'  definition turns out to be
\begin{defi} 
$S=\{P\in{\rm S}^3({k^2}^*):\exists
T_1,T_2\in T
\text{  s.t }
P=T_1+T_2 \text{ with  }T_1 ,T_2
 \text{  coprime}\}$.
\end{defi}
 Specializing to the space of binary cubics  a general theorem valid for   the 
symplectic covariants of the ${\mathfrak g_1}$ of any Heisenberg graded Lie algebra 
${\mathfrak g}$, we get the
\begin{theo} \label{sumofcubes}
 (i) Let $P\in S$ and let  $P=T_1+T_2$ with $T_1,T_2\in T$  coprime. Then $T_1,T_2$ 
are unique up to permutation. 
 
\medskip
\noindent(ii) Let  $P=T_1+T_2$ with $T_1,T_2\in T$. Then 
\begin{equation} Q_n(P)=\omega(T_1,T_2)^2.
\end{equation}
\medskip
\noindent(iii) Let $P\in {\cal O}_{[1]}$ and suppose $Q_n(P)=q^2$ with $q\in{k}^*$. 
Then
$$
T_1=\frac{1}{2}(P+\frac{1}{3q}\Psi(P)),\quad
T_2=\frac{1}{2}(P-\frac{1}{ 3q}\Psi(P))
$$
are coprime elements of $T$ such that $P=T_1+T_2$.
\end{theo}
\begin{demo}
For $k$ algebraically closed an argument that $P$ is a sum of cubes can be found in 
\cite[17-18]{Di}. The fact that $Q_n(P)=\omega(T_1,T_2)^2$ as well as (i) and (iii) are proved 
for general $k$ and for Heisenberg graded Lie algebras in \cite{Sl-St}.
\end{demo}
\begin{cor}
$S={\cal O}_{[1]}$.
\end{cor}
\begin{rema}
There is a natural bi-Lagrangian foliation of  ${\cal O}_{[1]}$ obtained by means of the decomposition $P = T_1 + T_2$.  Modulo some technicalities, if one fixes $T_2$ and varies over $T $ such that $\omega (T,T_2) = \omega(T_1,T_2)$ mod ${k^*}^2$, then does the same with $T_1$, one obtains a pair of foliations that are transverse and Lagrangian, for details see \cite{Sl-St}.
\end{rema}

Recall that elements of $T$ are, up to a scalar factor, cubes of linear forms. Hence a 
binary cubic $P$ is in $S$ iff
there exist a basis $\{\phi_1,\phi_2\}$ of ${k^2}^*$ and $\lambda_1,\lambda_2\in k^*$
 such that 
\begin{equation}\label{sumofcubes1}
P=\lambda_1\phi_1^3+\lambda_2\phi_2^3.
\end{equation}
The $\lambda_i$ and $\phi_i$ in this equation are not unique but the direct sum
decomposition 
$$
{k^2}^*=<\phi_1>\oplus<\phi_2>
$$ is canonically associated to $P$ as is described in the next result.
\begin{cor}  \label{eigvaluematchofmu}
(i) $P\in{\cal O}_{[1]}$ iff $\mu(P)\not=0$ is  diagonalisable over $k$, hence 
$\mu(P)$ is contained in a semisimple orbit.

\medskip
\noindent(ii) Let $P\in {\cal O}_{[1]}$ and let $\{\phi_1,\phi_2\}$ be a basis of 
${k^2}^*$. The following are equivalent:

\medskip
(a) There exist $\lambda_1,\lambda_2\in k^*$ such that 
$P=\lambda_1\phi_1^3+\lambda_2\phi_2^3$.

\medskip
(b) $\{\phi_1,\phi_2\}$ is a basis of eigenvectors of $\mu(P)$.

\bigskip
\noindent(iii) Let $P\in{\cal  O}_{[1]}$ and suppose $P=\lambda_1\phi_1^3+\lambda_2\phi_2^3$
where $\lambda_1,\lambda_2\in k^*$ and  $\{\phi_1,\phi_2\}$ is a basis of $k^2{^*}$. Then
if $q$ is the square root $\lambda_1\lambda_2\Omega(\phi_1,\phi_2)^3$ of $Q_n(P)$,
$$
\begin{array}{l}
\mu(P)\cdot \phi_1=-q\phi_1,\\
\mu(P)\cdot \phi_2=q\phi_2.
\end{array}
$$
\end{cor}
\begin{demo} (i): By Cayley-Hamilton and equation \eqref{normquart}, 
$$
0=\mu(P)^2+{\rm det}\mu(P)Id=\mu(P)^2-Q_n(P)Id.
$$
Hence $\mu(P)$ is diagonalisable over $k$ iff $Q_n(P)$ is a square in $k$.

\noindent(ii): Since there exists
$g\in \textrm{Sl}(2,k)$ with $<g\cdot \phi_1>=<x>$ and $<g\cdot \phi_2>=<y>$, we 
can assume without loss of generality that $\phi_1=x$ and $\phi_2=y$. Setting 
$P=ax^3+3bx^2y+3cxy^2+dy^3$, we have:
$\{x^3,y^3\}$ is a basis of eigenvectors of $\mu(P)$ iff $\mu(P)$ is diagonal iff (by  
equation \eqref{explicitmoment})
$$
bd-c^2=b^2-ac=0.
$$
This equation implies $b(ad-bc)=0$ and hence, since $Q_n(P)\not=0$, that $b=0$ and
$c^2=bd=0$. It follows that  $\{x^3,y^3\}$ is a basis of eigenvectors of $\mu(P)$ iff 
$b=c=0$ iff
$P=ax^3+dy^3$.

\noindent(iii):  As above, we can suppose without loss of generality that $P=ax^3+dy^3$ and then
$$
\mu(P)=
\begin{pmatrix}ad&0\\
0&-ad
\end{pmatrix}
$$
which implies $\mu(P)\cdot x=-adx$ and $\mu(P)\cdot y=ady$. This proves (iii) since $\Omega(x,y)=1$.
\end{demo}
\begin{cor}\label{fibmusquares}
(Fibres of $\mu$ on ${\cal O}_{[1]}$). Let $X\in\mathfrak{sl}(2,k)$ be diagonalisable over $k$, let $\pm q$ be its eigenvalues and let $\phi_+$ and $\phi_{-}$ be corresponding eigenvectors in $k^2{^*}$. Then
$$
\mu^{-1}(X)=\{
a\phi_{-}^3+{q\over a\Omega(\phi_{-},\phi_{+})^3}\,\phi_{+}^3:\,a\in k^*
\}.
$$
\end{cor}
\begin{demo} This follows from Corollary \ref{eigvaluematchofmu}(ii) and (iii).
\end{demo}
\vskip0,1cm

\centerline{{\bf Orbit parameters for ${\cal O}_{[1]}$ }}
\vskip0,1cm
For generic $k$ there will be many $\mathrm{Sl}(2,k)$ orbits on ${\cal O}_{[1]}$. So 
the first task is to obtain parameters for the orbits. For this the symplectic result 
Theorem  \ref{sumofcubes} leads to a new and effective method.
Let $P\in{\cal O}_{[1]}$.  Then as we have seen, there exist a unique {\it unordered} 
pair of elements $T_1,T_2$ in $T$ such that
\begin{align}
\nonumber
P&=T_1+T_2,\\
Q_n(P)&=\omega(T_1,T_2)^2.
\end{align}
Hence the map $I_{{\cal O}_{[1]}}:{\cal O}_{[1]}\rightarrow k^*\times_{Z_2} k^*/{k^*}^3$
\begin{equation}\label{diaginv}
I_{{\cal O}_{[1]}}(P)=[\omega(T_1,T_2),I_T(T_1)I_T(T_2)^{-1}]
\end{equation}
is well-defined where $k^*\times_{Z_2} k^*/{k^*}^3$ denotes the quotient of $k^*\times k^*/{k^*}^3$
by the $Z_2$-action
$$
-1\cdot(\lambda,\alpha)=(-\lambda,\alpha^{-1}).
$$
\begin{rema} The invariant $I_{{\cal O}_{[1]}}(\cdot)$ is symplectic not algebraic 
since its definition requires the symplectic form. We have not found this 
invariant for binary cubics in the literature.
\end{rema}
\begin{theo}\label{spofinv}
Let  $I_{{\cal O}_{[1]}}:{\cal O}_{[1]}\rightarrow k^*\times_{Z_2} k^*/{k^*}^3$ be
defined by  \eqref{diaginv} above.
\vskip0,1cm
\noindent(i) Let $P,P'\in {\cal O}_{[1]}$. Then 
$$ 
\textrm{Sl}(2,k)\cdot P'= \textrm{Sl}(2,k)\cdot P \iff
{{I_{{\cal O}_{[1]}}}}(P')={{I_{{\cal O}_{[1]}}}}(P).
$$

\noindent(ii) The map ${{I_{{\cal O}_{[1]}}}}$ induces a bijection  
$$
{\cal O}_{[1]}/\textrm{Sl}(2,k)\longleftrightarrow
k^*\times_{Z_2} k^*/{k^*}^3.$$

\noindent(iii) Let $P\in{\cal  O}_{[1]}$ and suppose $P=\lambda_1\phi_1^3+\lambda_2\phi_2^3$ where $\lambda_1,\lambda_2\in k^*$ and  $\{\phi_1,\phi_2\}$ is a basis of $k^2{^*}$. Let
$G_P=\{g\in\textrm{Sl}(2,k):\,g\cdot P=P\}$. Then
$$
G_P=\{g\in\textrm{Sl}(2,k):\,
\exists \mu\in k^*\text{ s.t. }g\cdot \phi_1=\mu\phi_1,g\cdot \phi_2={1\over\mu}\phi_2
\text{ and }\mu^3=1\}.
$$
\end{theo}
\begin{demo} (i): Since $\omega$ and $I_T$ are $\textrm{Sl}(2,k)$-invariant, it is clear from
\eqref{diaginv} that the map $I_{{\cal O}_{[1]}}:{\cal O}_{[1]}\rightarrow k^*\times_{Z_2} k^*/{k^*}^3$
factors  through the action of $\textrm{Sl}(2,k)$. To show that the induced map on orbit space is injective, suppose that $P$ and $P'$ are binary cubics such that $I_{{\cal O}_{[1]}}(P')=I_{{\cal O}_{[1]}}(P)$. First choose 
$g,g'\in \textrm{Sl}(2,k)$ such that 
\begin{align}
\nonumber
g\cdot P&=ax^3+by^3,\\
g'\cdot P'&=a'x^3+b'y^3.
\end{align}
From equations \eqref{bincubsymp} and \eqref{normquart} we have
$$
\omega(x^3,y^3)=1,\quad Q_n(P)=a^2b^2,\quad Q_n(P')={a'}^2{b'}^2.
$$
Hence $I_{{\cal O}_{[1]}}(P')=I_{{\cal O}_{[1]}}(P)$ implies
$$
[ab\,,[a\,][b]^{-1}\,]=[a'b'\,,[a'][b']^{-1}\,]
$$
in $k^*\times_{Z_2} k^*/{k^*}^3$. There are two possibilities:
\begin{itemize}
\item
$ab=a'b'$, \quad$[a\,][b]^{-1}=[a'][b']^{-1}$;
\item
$ab=-a'b'$, \quad$[a\,][b]^{-1}=[b'][a']^{-1}$.
\end{itemize}

In the first case, we have
$$
[ab][a\,][b]^{-1}=[a'b'][a'][b']^{-1},
$$
hence $[a^2]=[{a'}^2]$ and so $[a]=[a']$ as the group $k^*/{k^*}^3$ is of exponent $3$.
Thus there exists $r\in k^*$  such that $a'=r^3a$ and 
$b'={1\over r^3}b$.  If we define $h\in GL(2,k)$ by
$$
h\cdot x=rx,\quad h\cdot y= {1\over r}y,
$$
it is clear that $h\in \textrm{Sl}(2,k)$ and $h\cdot(g\cdot P)=g'\cdot P'$. Hence $P$ 
and $P'$ are in the same $\textrm{Sl}(2,k)$-orbit.

In the second case, we have $[a^2]=[{b'}^2]$, $[a]=[b']$ and there exists $r\in k^*$ such that
$b'=r^3a$ and $a'=-{1\over r^3}b$. If we define $h\in GL(2,k)$ by
$$
h\cdot x=ry,\quad h\cdot y= -{1\over r}x,
$$
it is clear that $h\in \textrm{Sl}(2,k)$ and $h\cdot(g\cdot P)=g'\cdot P'$.  Hence $P$ 
and $P'$ are in the same $\textrm{Sl}(2,k)$-orbit
and we have proved that $I_{{\cal O}_{[1]}}:{\cal O}_{[1]}\rightarrow k^*\times_{Z_2} k^*/{k^*}^3$ separates $\textrm{Sl}(2,k)$-orbits.

To prove (ii) it remains to prove that $I_{{\cal O}_{[1]}}:{\cal O}_{[1]}\rightarrow k^*\times_{Z_2} k^*/{k^*}^3$ is surjective.  Let $[q,[\alpha]]\in k^*\times_{Z_2} k^*/{k^*}^3$ and consider the binary cubic
$$
P={1\over q\alpha}x^3+q^2\alpha y^3.
$$
Then
$$
I_{{\cal O}_{[1]}}(P)=
[q,\,[{1\over q\alpha}]\,[q^2\alpha\,]^{-1}]
=[q,[{1\over q^3\alpha^2}]]=[q,[\alpha]].
$$
 and so $I_{{\cal O}_{[1]}}:{\cal 
O}_{[1]}\rightarrow k^*\times_{Z_2} k^*/{k^*}^3$ is
 surjective. This completes the proof of (ii).
 
 To prove (iii), recall that the representation  $P = \lambda_1 \phi_1^3 + \lambda_2 \phi_2^3$ is unique up to permutation. Then $g\cdot P = P$  leads to two cases: 
 \begin{itemize}
\item
 $g\cdot(\lambda_1 \phi_1^3)=\lambda_1 \phi_1^3$ and $g\cdot(\lambda_2 \phi_2^3)=\lambda_2 \phi_2^3$;
\item
$g\cdot(\lambda_1 \phi_1^3)=\lambda_2 \phi_2^3$ and $g\cdot(\lambda_2 \phi_2^3)=\lambda_1 \phi_1^3$.
\end{itemize}
In the first case, $g\cdot\phi_i = j_i\phi_i$ where $j_i^3 = 1$ and since $g \in \textrm{Sl}(2,k)$, we must have $j_1j_2=1$. In the second case, there exist $r,s\in k^*$ such that $g\cdot \phi_1 = r\phi_2$, $g\cdot \phi_2 =s \phi_1$, $\lambda_1r^3=\lambda_2$, $\lambda_2s^3=\lambda_1$ and $rs=-1$. Hence $(rs)^3=1$ and $rs=-1$ which is impossible and this case does not occur.
 \end{demo}

\vskip0,1cm
\centerline{{\bf Properties of orbit space}}
\vskip0,1cm
We will use the parametrisation
$$
I_{{\cal O}_{[1]}}:{{\cal O}_{[1]}}/\textrm{Sl}(2,k)\longleftrightarrow  k^*\times_{Z_2}k^*/{k^*}^3.
$$
 to study orbit space. The parameter space has two
natural maps
\begin{equation} sq:k^*\times_{Z_2}k^*/{k^*}^3\rightarrow{k^*}^2,\quad 
sq([q,\alpha])=q^2,
\end{equation}
and
\begin{equation} t:k^*\times_{Z_2}k^*/{k^*}^3\rightarrow(k^*/{k^*}^3)/Z_2,\quad 
t([q,\alpha])=[\alpha]
\end{equation}
corresponding to projection onto the orbit spaces of the two factors.
 We then have the following diagram:
\begin{equation}\label{nonunitarydiag}
 \xymatrix
{
&k^*\times_{Z_2}k^*/{k^*}^3\ar[dl]^{sq}\ar[dr]^{t}&\\
{k^*}^2\ar@/^/[ur]^{e}&&(k^*/{k^*}^3)/Z_2.\\
}
 \end{equation}

 \vskip0,1cm
The map
\begin{equation}\label{nonunitarysq}
sq:k^*\times_{ Z_2}k^*/{k^*}^3
\rightarrow {k^*}^2
\end{equation}
is the fibration associated to the principal $Z_2$-fibration
$$
{k}^*\rightarrow{k^*}^2
$$
and the action of $Z_2$ on $k^*/{k^*}^3$ by inversion. Since $Z_2$ acts by automorphisms, the fibre $sq^{-1}(q^2)$ over any point  $q^2\in {k^*}^2$ has
a natural group structure
\begin{equation}\label{gpstrnonunitary}
[q,\alpha]\times[q,\beta]=[q,\alpha\beta]
\end{equation}
independent of the choice of square root $q$ of $q^2$. Taking the identity at each point, we get a canonical section 
$e:{k^*}^2\rightarrow k^*\times_{ Z_2}k^*/{k^*}^3$ of 
\eqref{nonunitarysq} given by
\begin{equation}\label{iddiscfixnonunitary}
e(q^2)=[q,1]
\end{equation}
but, although each fibre is a group isomorphic to $k^*/{k^*}^3$, the fibration
\eqref{nonunitarysq} is not in general isomorphic to the product 
$$
{k^*}^2\times k^*/{k^*}^3
\rightarrow{k^*}^2.
$$

To translate the above features of orbit space into  more concrete statements about binary cubics over $k$, note that the map
$sq$ is essentially the quartic $Q_n$ since for all $P\in{\cal O}_{[1]}$,
$$
sq({{I_{{\cal O}_{[1]}}}}(P))=Q_n(P).
$$
\begin{theo}\label{fixeddiscnonunitary}
Let $M\in {k^*}^2$, let 
$$
{\cal O}_M=\{P\in S^3({k^2}^*):\,Q_n(P)=M\}
$$
and let ${\cal O}_M/\textrm{Sl}(2,k)$ be the  space of ${Sl}(2,k)$-orbits in ${\cal O}_M$. 

\medskip
\noindent(i) The map ${I_{{\cal O}_{[1]}}}:{\cal O}_{[1]}
\rightarrow k^*\times_{ Z_2}k^*/{k^*}^3$ induces a bijection
$$
{\cal O}_M/\textrm{Sl}(2,k)\longleftrightarrow
{{sq}}^{-1}(M)
$$
and, by pullback of \eqref{gpstrnonunitary}, a group structure on ${\cal O}_M/\textrm{Sl}(2,k)$.

\medskip
\noindent(ii) As groups, ${\cal O}_M/\textrm{Sl}(2,k)\cong k^*/{k^*}^3$.

\medskip
\noindent(iii) Let $q\in k^*$ be a square root of $M$. The identity element of ${\cal O}_M/\textrm{Sl}(2,k)$ is characterised by:
$$
\textrm{Sl}(2,k)\cdot P=1\Leftrightarrow P\text{ is reducible over }k
\Leftrightarrow
{I_{{\cal O}_{[1]}}}(P)=[q,1].
$$
\end{theo} 
\begin{demo} Parts (i) and (ii) follow from the discussion above. Part (iii) follows
from Theorem \ref{reduc}(i) and equation \eqref{iddiscfixnonunitary}.
\end{demo}
\begin{rema} \label{uniqredorb1}
From the Corollary it follows that if the classical discriminant is a nonzero square there is a unique $\textrm{Sl}(2,k)$ orbit consisting of reducible polynomials. We remove the `square' restriction  
in Corollary \ref{uniqredorb2}. In particular, over an algebraically closed field there is only one orbit of fixed nonzero discriminant. 
\end{rema}

To finish this section we briefly discuss the map
$t:k^*\times_{Z_2}k^*/{k^*}^3\rightarrow (k^*/{k^*}^3)/Z_2$ in diagram \eqref{nonunitarydiag}
given by
$$
t([q,\alpha])=[\alpha].
$$
This a fibration with fibre $k^*$ outside the identity coset $[1]$ but  
 $$
t^{-1}([1])=e({k^*}^2)
 $$
  is a `singular fibre'.
   There is a $k^*$-action:
\begin{equation}\label{spin}
\lambda\cdot [q,\alpha]=[\lambda q,\alpha]
\end{equation}
which maps fibres of $sq$ to fibres of $sq$:
$$
{{sq}}([q',\alpha'])={{sq}}([q,\alpha])\Rightarrow {{sq}}(\lambda\cdot[q',\alpha'])={{sq}}(\lambda\cdot[q,\alpha]),
$$
and whose orbits are exactly the fibres of $t$:
$$
t([q',\alpha'])=t([q,\alpha])\Leftrightarrow \exists\lambda\in k^*\text{ s.t. }
[q',\alpha']=\lambda\cdot [q,\alpha].
$$
Isotropy for this action is given by: 
$
Isot_{k^*}([q,\alpha])=\left\{
\begin{array}{lr}
1&\text{if }\alpha\not=1\\
\{\pm 1\}&\text{if }\alpha=1.
\end{array}
\right.
$

It would be interesting to interpret these features of orbit space in terms of the original binary cubics. Conversely, one can also identify actions on the orbits in terms of their orbit parameters.
For example, the commutant of $\textrm{Sl}(2,k)$ in $\textrm{Gl}({\rm S}^3({k^2}^*))$  acts on orbit space.  This gives the action  
$$
\lambda\cdot'[q,\alpha]=[\lambda^2 q,\alpha]
$$
of $k^*$ on $k^*\times_{Z_2}k^*/{k^*}^3$ which is the square of the action \eqref{spin}.
Another example is obtained from  $\Psi:{\rm S}^3({k^2}^*)\rightarrow {\rm S}^3({k^2}^*)$ which, since it commutes with the 
action of
$\textrm{Sl}(2,k)$, induces a map from $k^*\times_{Z_2}k^*/{k^*}^3$ to itself. 
This is easily seen to be given by
\begin{equation}
[q,\alpha]\mapsto [-q^3,[q]\alpha],
\end{equation}
where $[q]$ denotes the class of $q$ in $k^*/{k^*}^3$.

\vskip0,1cm
\centerline{{\bf Reducibility and factorisation}}
\vskip0,1cm
\begin{theo} \label{reduc}
Let $P\in S$ and let $\{\phi_1,\phi_2\}$ be a basis of ${k^2}^*$ such that  
$P=\lambda_1\phi_1^3+\lambda_2\phi_2^3$
with
$\lambda_1,\lambda_2\in k^*$.
Let $q\in k^*$ be a square root of $Q_n(P)$. The following are equivalent:

\medskip
\noindent(a) $P$ is reducible over $k$.

\medskip
\noindent(b) $\frac{\lambda_1}{ \lambda_2}$ is a cube in $k^*$.

\medskip
\noindent(c) $q\lambda_1$ is a cube in $k^*$.

\medskip
\noindent(d) $q\lambda_2$ is a cube in $k^*$.

\medskip
\noindent(e) There is a basis $\{\phi'_1,\phi'_2\}$ of ${k^2}^*$ such that 
$P=\frac{1}{q}({\phi'}_1^3+{\phi'}_2^3)$.
\end{theo}
\begin{demo} \noindent$(a)\Rightarrow(b)$: Suppose $P$ is reducible over $k$. Then for 
all $g\in \textrm{Sl}(2,k)$,
$$
g\cdot P=\lambda_1(g\cdot\phi_1)^3+\lambda_2(g\cdot\phi_2)^3
$$
is also reducible over $k$.  Since $\phi_1,\phi_2$ form a basis of ${k^2}^*$, we can 
choose $g$ such that 
$g\cdot\phi_1=x$ and $g\cdot\phi_2=\rho y$ for some $\rho\in k^*$ so that 
$$
\lambda_1x^3+\lambda_2\rho^3y^3
$$
is reducible over $k$. Hence there exist $a,b,c,d,e\in k$ such that
$$
\lambda_1x^3+\lambda_2\rho^3y^3=(ax+by)(cx^2+dxy+ey^2)
$$
which gives the system
$$
\begin{array}{rl}
\lambda_1=ac,&0=ad+bc,\\
\lambda_2\rho^3=be,&0=ae+bd.
\end{array}
$$
Since $\lambda_1$ and $\lambda_2$ are nonzero, it follows that $a,b,c,d,e$ are nonzero
and, since $c=-\frac{ad}{b}$ and $e=-\frac{bd}{a}$  we get $\frac{\lambda_1}{ 
\lambda_2}=(\rho\frac{a}{b})^3$.

\noindent$(b)\Rightarrow(a)$: Suppose $\frac{\lambda_1}{\lambda_2}=r^3$ with $r\in 
k^*$. Then
\begin{equation}\label{firstfact}
P=\lambda_2(r^3\phi_1^3+\phi_2^3)=\lambda_2(r\phi_1+\phi_2)(r^2\phi_1^2+r\phi_1\phi_2+
\phi_2^2)
\end{equation}
and $P$ is reducible over $k$. 

\noindent$(b)\Leftrightarrow(c)\Leftrightarrow(d)$:  Set $\nu_1=q\lambda_1$ and 
$\nu_2=q\lambda_2$. By Proposition \ref{constraint}, there exists $s\in k^*$ such that 
$\nu_1\nu_2=s^3$. Hence if any one of the three numbers  $\nu_1, \nu_2, 
\frac{\nu_1}{\nu_2}= \frac{\lambda_1}{\lambda_2}$ is a cube so are the other two since 
formally
$$
\nu_1=\left(\frac{\nu_1}{\sqrt[3]{\frac{\nu_1}{\nu_2}}\sqrt[3]{\nu_1\nu_2}}\right)^3,\quad
\nu_2=\left(\sqrt[3]{\frac{\nu_1}{\nu_2} }  \frac 
{\nu_2}{\sqrt[3]{\nu_1\nu_2}}\right)^3,\quad
\nu_2=\left(\frac{ \sqrt[3]{\nu_1\nu_2}}{\sqrt[3]{\nu_1}} \right)^3.
$$

\noindent$(a)\Rightarrow(e)$: If $P$ is reducible we have just proved that
 there exists $r\in k^*$  and $s\in k^*$such that
$\lambda_1=\frac{1}{q}r^3$ and $\lambda_2=\frac{1}{q}s^3$.
Set   $\phi'_1=r\phi_1$ and $\phi'_2=s\phi_2$. Then
$$
P=\lambda_1\phi_1^3+\lambda_2\phi_2^3=\frac{1}{q}\left( {\phi'_1}^3+{\phi'_2}^3 
\right)
$$
which proves $(e)$.

\noindent$(e)\Rightarrow(a)$: Evident since ${\phi'}_1+{\phi'}_2$ divides 
${\phi'}_1^3+{\phi'}_2^3$.
\end{demo}

\begin{cor} \label{factors}
Let $P\in S$ be reducible and let $\{\phi'_1,\phi'_2\}$ be a  basis of ${k^2}^*$ such 
that 
 $P=\frac{1}{q}({\phi'}_1^3+{\phi'}_2^3)$.

\medskip
\noindent(a) If   $-3$ is not a square in $k$, then
$$
P=\frac{1}{q}(\phi'_1+\phi'_2){(\phi'_1}^2-\phi'_1\phi'_2+{\phi'_2}^2)
$$
and ${\phi'_1}^2-\phi'_1\phi'_2+{\phi'_2}^2$ is irreducible over $k$.

\medskip
\noindent(b) If  $-3$ is  a square in $k$, then
\begin{equation}
P=\frac{1}{q}(\phi'_1+\phi'_2)(j\phi'_1+j^{-1}\phi'_2)(j^2\phi'_1+j^{-2}\phi'_2)
\end{equation}
where $j=\frac{1}{2}(-1+\sqrt{-3})$. The factors of $P$ are pairwise independent.
\end{cor}
To a certain extent, we can normalise bases of  ${k^2}^*$ satisfying  Theorem 
\ref{sumofcubes}(e).
\begin{cor}\label{normedsum}
Let $P\in S$. 

\medskip
\noindent(a)  $P$ is reducible iff there is a basis $\{\phi'_1,\phi'_2\}$ of ${k^2}^*$ 
such that $P=\frac{1}{q}({\phi'}_1^3+{\phi'}_2^3)$ and
$\Omega({\phi'}_1,{\phi'}_2)=q $.

\medskip
\noindent(b If  $\{\phi'_1,\phi'_2\}$ and  $\{\phi''_1,\phi''_2\}$
 are two    bases of ${k^2}^*$ satisfying (a),
   there exists a cube root of unity $j\in k^*$ such that ${\phi''_1}=j{\phi'_1}$ and 
${\phi''_2}=j^{-1}{\phi'_2}$.
\end{cor}
\begin{demo}
Choose a basis  $\{\phi_1,\phi_2\}$ of ${k^2}^*$ and $\lambda_1,\lambda_2\in k^*$ such 
that $P=\lambda_1\phi_1^3+\lambda_2\phi_2^3$ and let $q=\lambda_1\lambda_2\Omega(\phi_1,\phi_2)^3$.  If $P$ is reducible, by Theorem 
\ref{sumofcubes}, there exists $r\in k^*$ such that
$\lambda_1=\frac{1}{q}r^3$.  Set $s=\frac{q}{r\Omega(\phi_1,\phi_2)}$,  
$\phi'_1=r\phi_1$
and $\phi'_2=s\phi_2$. Then $\Omega({\phi'}_1,{\phi'}_2)=q$ and
$$
s^3=\left(\frac{q}{r\Omega(\phi_1,\phi_2)}\right)^3
=\frac{1}{r^3}\left(\frac{q}{\Omega(\phi_1,\phi_2)}\right)^3
=\frac{1}{q\lambda_1}\left(q\lambda_1\right)\left(q\lambda_2\right)=
q\lambda_2.
$$
Hence
$$
P=\lambda_1\phi_1^3+\lambda_2\phi_2^3=\frac{1}{q}\left( {\phi'_1}^3+{\phi'_2}^3 
\right).
$$
In the classical literature on cubics this is called the Vi\`ete Substitution.

Conversely,  if there is a basis $\{\phi'_1,\phi'_2\}$ of ${k^2}^*$ such that 
$P=\frac{1}{q}({\phi'}_1^3+{\phi'}_2^3)$,
then $\phi'_1+\phi'_2$ divides $P$ and $P$ is reducible.

 To prove (b), note first that by Theorem \ref{sumofcubes}(a), we have
either ${\phi''_1}^3={\phi'_1}^3$ and ${\phi''_2}^3={\phi'}_2^3$ or 
${\phi''_1}^3={\phi'_2}^3$ and ${\phi''_2}^3={\phi'}_1^3$. 

In the first case,
by unique factorisation, there
exist  cube roots of unity $j_1,j_2$ such that  ${\phi''_1}=j_1{\phi'_1}$, 
${\phi''_2}=j_2{\phi'_2}$ and $j_1j_2=1$.
This is exactly what we want to prove.

In the second case, there
exist  cube roots of unity $j_1,j_2$ such that  ${\phi''_1}=j_1{\phi'_2}$, 
${\phi''_2}=j_2{\phi'_1}$ and $j_1j_2=-1$. This is impossible since $(j_1j_2)^3=1$.
\end{demo}
\vskip0,1cm
\centerline{{\bf Explicit formulae for  $I_{{\cal O}_{[1]}}$ and Cardano-Tartaglia 
formulae}} 
\vskip0,1cm
\begin{prop} \label{explicitI}
Let $P=ax^3+3bx^2y+3cxy^2+dy^3$ be an element of ${\cal O}_{[1]}$,   let
$q\in k^*$ be a square root of $Q_n(P)$ and define $\alpha,\beta,\gamma$ and $\delta$ 
in $k$ by
$$
\mu(P)=
\begin{pmatrix}
(ad-bc)&2(bd-c^2)\\
2(b^2-ac)&-(ad-bc)
\end{pmatrix}
=
\begin{pmatrix}
\alpha&\beta\\
\gamma&\delta
\end{pmatrix}.
$$
Then $P=\lambda_1\phi_1^3+\lambda_2\phi_2^3$ and
$I_{{\cal O}_{[1]}}(P)=[\,q\,,[\lambda_1][\lambda_2]^{-1}\,]$ where:

\medskip
(i) If $\beta=\gamma=0$,
$$
\begin{array}{rl}
\lambda_1=a,&
\phi_1=x,\\
\lambda_2=d,&
\phi_2=y,
\end{array}
\quad ad= q.
$$

(ii) If $\gamma\not=0$,
$$
\begin{array}{ll}
\lambda_1=\frac{1}{2q}(\alpha+q)a+\frac{\gamma}{2q}b,&
\phi_1=x-(\frac{\alpha-q}{\gamma})y,\\
\lambda_2=-\frac{1}{2q}(\alpha-q)a-\frac{\gamma}{2q}b,&
\phi_2=x-(\frac{\alpha+q}{\gamma})y,
\end{array}
\quad\Omega(\phi_1,\phi_2)=-\frac{2q}{\gamma}.
$$

(iii)  If $\beta\not=0$,
$$
\begin{array}{ll}
\lambda_1=\frac{\beta}{ 2q}c-\frac{1}{2q}(\alpha-q)d,&
\phi_1=(\frac{\alpha+q}{\beta})x+y,\\
\lambda_2=-\frac{\beta}{2q}c+\frac{1}{2q}(\alpha+q)d,&
\phi_2=(\frac{\alpha-q}{\beta})x+y,
\end{array}
\quad\Omega(\phi_1,\phi_2)=\frac{2q}{\beta}.
$$
\end{prop}
If $P\in{\cal O}_{[1]}$ is reducible we can use these formulae together with Theorem 
\ref{reduc} and Corollary \ref{factors} to get an explicit formula for a linear factor 
of $P$ in terms of the coefficients of $P$, a square root $q$ of
$Q_n(P)$ and a cube root $r$ of $q\lambda_1$. Recall that the existence of 
a cube root  of $q\lambda_1$ in $k$ is a necessary and sufficient condition 
for $P$ to be reducible over $k$.
\begin{prop} \label{CardanoTartaglia}
Let $P=ax^3+3bx^2y+3cxy^2+dy^3\in {\cal O}_{[1]}$ be reducible,   let
$q\in k^*$ be a square root of $Q_n(P)$ and suppose $ad\not=0$. 

\medskip
(i) If $\beta=\gamma=0$, let $r$ be a cube root of $qa$ and let 
$s=\frac{q}{r}.$ Then
$$
rx+sy
$$

divides $P$.

(ii) If $\gamma\not=0$, let $r$ be a cube root of $(\alpha+q)a+\gamma b$ and 
let $s=-\frac{\gamma}{r}$.
Then
$$
x+\left(\frac{r-s+b}{a}\right)y
$$

divides $P$.

(iii)  If $\beta\not=0$, let $r$ be a cube root of $\beta c-(\alpha-q)d$ and 
let $s=\frac{\beta}{ r}$. Then
$$
\left(\frac{s-r+c}{d}\right)x+y
$$

divides $P$.

\end{prop}
\begin{demo}
Since $P$ is reducible,  there exists a basis $\phi'_1,\phi'_2$ of ${k^2}^*$ such that
$P=\frac{1}{q}({\phi'}_1^3+{\phi'}_2^3)$ (cf Theorem \ref{reduc} ) and then 
${\phi'}_1+{\phi'}_2$ divides $P$.  As shown in the proof
of Corollary \ref{normedsum}$(a)$, we can take ${\phi'}_1=r\phi_1$ and 
${\phi'}_2=s\phi_2$ where  $r$ is a cube root
of $q\lambda_1$, $s=\frac{q}{r\Omega(\phi_1,\phi_2)}$ and
$\phi_1,\phi_2,\lambda_1$ are given by Proposition \ref{explicitI}.
The explicit formulae in the three cases are:

\medskip
(a) $\beta=\gamma=0$: \quad $r$ is a cube root of $qa$, $rs=q$ 
and
$
\phi'_1=rx,\,
\phi'_2=sy;
$

\medskip
(b) $\gamma\not=0$:\quad$r$ is a cube root of $\frac{(\alpha+q)a+\gamma b}{2}$, 
$s=-\frac{\gamma}{2 r}$ and
$$
\phi'_1=rx+\frac{1}{2s}(\alpha-q)y,\quad
\phi'_2=sx+\frac{1}{2r}(\alpha+q)y;
$$

(c)   $\beta\not=0$:\quad$r$ is a cube root of $\frac{\beta c-(\alpha-q)d}{2}$, 
$s=\frac{\beta}{ 2r}$ and
$$
\phi'_1=\frac{1}{2s}(\alpha+q)x+ry,\quad
\phi'_2=\frac{1}{2r}(\alpha-q)x+sy.
$$
Calculating ${\phi'}_1+{\phi'}_2$ in the first case obviously gives (i). In the second 
case we have
\begin{align}
\nonumber
{\phi'}_1+{\phi'}_2&=(r+s)x+\left(\frac{1}{2s}(\alpha-q)+\frac{1}{ 
2r}(\alpha+q)\right)y\\
&=(r+s)x+\left(\frac{1}{2sa}(-2s^3-\gamma b)+\frac{1}{2ra}(2r^3-\gamma b)\right)y
\end{align}
since $r^3=q\lambda_1=\frac{(\alpha+q)a+\gamma b}{2}$ and 
$s^3=q\lambda_2=\frac{-(\alpha-q)a-\gamma b}{2}$. Simplifying the 
coefficient of $y$ we get
$$
\frac{1}{ 2sa}(-2s^3-\gamma b)+\frac{1}{ 2ra}(2r^3-\gamma b)=
\frac{1}{ a}\left(r^2-s^2-\frac{b\gamma}{2}(\frac{1}{ r}+\frac{1}{ s})
\right)
=(r+s)\frac{r-s+b}{ a}
$$
since $2rs=-\gamma$, and this implies (ii). Similarly, (iii) follows from (c).
\end{demo}

As an application of the above results, consider the homogeneous Cardano-Tartaglia 
polynomial
$$
P=x^3+pxy^2+qy^3
$$
over a field $k$ of characteristic not $2$ or $3$.
Assume $p\not=0$ and $q\not=0$ so that factorising $P$ is  a nontrivial problem.
Then 
$$
\mu(P)=\begin{pmatrix}
q&-2\frac{p^2}{ 9}\\
-2\frac{p}{ 3}&-q
\end{pmatrix},\quad Q_n(P)=(q^2+4\frac{p^3}{27}).
$$
To be able to apply our approach we assume $Q_n(P)$ has a square root  in $ k^*$ which 
we denote $\sqrt{q^2+4\frac{p^3}{27}}$. 
Then by Theorem \ref{reduc} and Proposition \ref{explicitI}(ii), $P$ is reducible iff
$$
\frac{q}{ 2}+\sqrt{\frac{q^2}{ 4}+\frac{p^3}{27}}\quad
\text{or}
\quad
-\frac{q}{ 2}+\sqrt{\frac{q^2}{ 4}+\frac{p^3}{ 27}}
$$
has a cube root in $k$.

If this is the case, then Proposition \ref{CardanoTartaglia} (ii)  implies that 
$x+(r-s)y$ divides $P$ where $r$ is a cube root of 
$
\frac{q}{ 2}+\sqrt{\frac{q^2}{ 4}+\frac{p^3}{ 27}}
$
and $s$ is the cube root $\frac{p}{ 3r}$ of $-\frac{q}{ 2}+\sqrt{\frac{q^2}{ 
4}+\frac{p^3}{ 27}}$. Hence, with the obvious notation,
$$
\frac{p}{3\left(\sqrt[3]{\frac{q}{ 2}+\sqrt{\frac{q^2}{ 4}+\frac{p^3}{27}}}\right)}-
\sqrt[3]{\frac{q}{ 2}+\sqrt{\frac{q^2}{ 4}+\frac{p^3}{ 27}}}
$$
 is a root of the inhomogeneous cubic $x^3+px^2+q$ and this is the classical   
Cardano-Tartaglia formula. If
$k={\bb R}$, this can be written
$$
s-r=\sqrt[3]{-\frac{q}{ 2}+\sqrt{\frac{q^2}{ 4}+\frac{p^3}{ 27}}}-\sqrt[3]{\frac{q}{ 
2}+\sqrt{\frac{q^2}{ 4}+\frac{p^3}{27}}}
$$
since cube roots are unique.

\subsection
{Symplectic covariants and sums of coprime cubes in quadratic extensions}

In this article  we have until now considered only binary cubics $P$ such that $Q_n(P)$ is a square in $k$. In this section we will study binary cubics $P$ such that $Q_n(P)$ is a square in a  fixed quadratic extension of $k$.

Let $\hat{k}$ be a quadratic extension of $k$. Recall that since char$(k)\not=2$, the extension $\hat{k}/k$ is Galois and  the Galois group ${\rm Gal}(\hat{k}/k)$ is isomorphic to $Z_2$. The Galois group ${\rm Gal}(\hat{k}/k)$ acts naturally on any  space over $\hat{k}$ obtained by base extension of a   space over $k$ and its fixed point set is the original space over $k$. We
always denote the action of the generator of  ${\rm Gal}(\hat{k}/k)$ by
$x\mapsto\bar{x}$ and we denote by $\hat{\Omega}$ and $\hat{\omega}$ respectively
the symplectic forms on ${{\hat{k}}^2}{^*}$ and
${S^3}({{\hat{k}}^2}{^*})$ obtained  by base extension of $\Omega$ and $\omega$. The quartic on ${S^3}({{\hat{k}}^2}{^*})$ obtained by base extension of $Q_n$ will be denoted $\widehat{Q_n}$ and we set
$$
\widehat{\cal O}_{[1]}=\{P\in{S^3}({{\hat{k}}^2}{^*}):\,\widehat{Q_n}(P)\in\hat{k}{^*}{^2}\}.
$$
Finally, 
let ${\rm Im}\,\hat{k}=\{\lambda\in \hat{k}:\,\bar{\lambda}=-\lambda\}$ and 
let $\widehat{T}\subseteq{S^3}({{\hat{k}}^2}{^*})$ be the set of nonzero binary cubics over $\hat{k}$ which have a triple root over $\hat{k}$. 
\begin{rema}
Note  that $({\rm Im}\,{\hat{k}}^*)^2\subseteq k^*$ is the inverse image under $k^*\rightarrow{k^*}^2$ of a single nontrivial square class in $k^*/{k^*}^2$. Conversely, a  nontrivial square class in $k^*/{k^*}^2$ determines up to isomorphism a quadratic extension of $k$ with this property.
\end{rema}
This notation out of the way, we make a symplectic definition
$$
{\cal O}({\hat{k}})=\{P\in S^3({k^2}^*):\,\hat{k}\text{ is a splitting field of } x^2-Q_n(P)\}
$$
and an algebraic definition
$$
S({\hat{k}})=\{P\in S^3({k^2}^*):\,\exists T\in \widehat{T}\text{ s.t. }P=T+\bar{T}\text{ with }T,\bar{T}\text{ coprime}\}.
$$
\begin{prop}
${\cal O}({\hat{k}})=S({\hat{k}})$.
\end{prop}
\begin{demo}
Let $P\in{\cal O}({\hat{k}})$. Then $Q_n(P)$ has two square roots in $\hat{k}$ but no square roots in $k$ since
$\hat{k}$  is a splitting field of  $x^2-Q_n(P)$. By Theorem
 \ref{sumofcubes}, there exists $T_1,T_2\in \widehat{T}$  such that $P=T_1+T_2$ and the square roots of $Q_n(P)$ are $\pm{\hat{\omega}}(T_1,T_2)$. Since $\bar{P}=P$ and since $T_1$ and $T_2$ are unique up to permutation, we have either $\bar{T_1}=T_1$ and $\bar{T_2}=T_2$ or $\bar{T_1}=T_2$ and $\bar{T_2}=T_1$. In the first case, 
$$
\overline{\hat{\omega}(T_1,T_2)}=\hat{\omega}(\bar{T_1},\bar{T_2})=\hat{\omega}(T_1,T_2),
$$
so $\hat{\omega}(T_1,T_2)\in k$ and  $Q_n(P)$ has a square root in  $k$ which is a contradiction. Hence $P=T_1+\bar{T_1}$. To prove that $T_1$ and $\bar{T_1}$ are coprime,
write $T_1=\lambda\alpha^3$ where $\lambda\in\hat{k}$ and $\alpha\in{\hat{k}^2}{^*}$. Then, by unique factorisation, 
$T_1$ and $\bar{T_1}$ are not coprime iff $\alpha$ and $\overline{\alpha}$ are proportional.
But then $\hat{\omega}(T_1,\bar{T}_1)=0$ and $Q_n(P)=0$ has a square root in  $k$. Hence 
$T_1$ and $\bar{T_1}$ are coprime and $P\in S({\hat{k}})$.

To prove inclusion in the opposite direction, suppose $P\in S({\hat{k}})$ and let $P=
T+\bar{T}$ with $T,\bar{T}$ coprime and $T\in\widehat{T}$. Note that $P\not=0$ since otherwise
$T$ and $\bar{T}$ would not be coprime.
By Theorem
 \ref{sumofcubes}, we have $Q_n(P)=(\hat{\omega}(T,\bar{T}))^2$ and $Q_n(P)$ has two square roots ${\hat{\omega}}(T,\bar{T})$ in
$\hat{k}$.  Let $T=\lambda\alpha^3$ where $\lambda\in\hat{k}^*$ and $\alpha\in{\hat{k}^2}{^*}$. 
As we saw above, $T$ and $\bar{T}$ are  coprime implies $\alpha$ and $\overline{\alpha}$ are not proportional, and this is equivalent to $\hat{\Omega}(\alpha,\overline{\alpha})\not=0$
since ${\rm dim}\,{\hat{k}^2}{^*}=2$.
From
$$
\hat{\omega}(T,\bar{T})=\lambda\bar{\lambda}(\hat{\Omega}(\alpha,\overline{\alpha}))^3
$$
it follows that $\hat{\omega}(T,\bar{T})\not=0$. On the other hand,
$$
\overline{\hat{\omega}(T,\bar{T})}=\hat{\omega}(\bar{T},T)=-\hat{\omega}(T,\bar{T})
$$
and $\hat{\omega}(T,\bar{T})$ is pure  imaginary. Hence
the square roots $\pm{\hat{\omega}}(T,\bar{T})$ of $Q_n(P)$are not in $k$ and $\hat{k}$ is a splitting field of $x^2-Q_n(P)$.
\end{demo}
\begin{prop}(Fibres of $\mu$ on ${\cal O}({\hat{k}})$). 
Let $X\in\mathfrak{sl}(2,k)$ be such that $-\textrm{det}\,X\in ({\rm Im}\,{\hat{k}}^*)^2$ and
$\nu_{\textrm{det}\,X}(X)=[2]$.  Let $q,\bar{q} \in {\rm Im}\,{\hat{k}}^*$
 be its eigenvalues and let $\phi$ and $\bar{\phi}$ be corresponding eigenvectors in $\hat{k}^2{^*}$. 
 
 \medskip
 \noindent(i) There exists $a\in\hat{k}^*$ such that $a\bar{a}\,\Omega(\bar{\phi},\phi)^3=q$.
 
 \medskip
 \noindent(ii)
 $$
\mu^{-1}(X)=\{
ua{\phi}^3+\bar{u}\bar{a}\,\bar{\phi}^3:\,u\in \hat{k}^*\text{ and }u\bar{u}=1
\}.
$$
\end{prop}
\begin{demo} 
Recall that $\nu_{\textrm{det}\,X}(X)=[2]$ is a necessary and sufficient condition for $X$ to be in the image of $\mu$ (cf Theorem \ref{Immom}). Since $\phi+\bar{\phi}$ is not
an eigenvector of $X$, we have
$$
[2]=[\Omega(\phi+\bar{\phi},X\cdot\phi+X\cdot\bar{\phi})]=[-2q\Omega(\phi,\bar{\phi})].
$$
Hence there exists $\alpha\in\hat{k}^*$ such that $\alpha\bar{\alpha}=q\Omega(\bar{\phi},\phi)$ and then $a={q^2\over \alpha^3}$ is  a solution of (i).

By Corollary \ref{fibmusquares}, the fibre of the $\hat{k}$-moment map $\hat{\mu}:S^3({\hat{k}^2}{^*})\rightarrow\textrm{sl}(2,\hat{k})$ is
$$
\hat{\mu}^{-1}(X)=\{
c\bar{\phi}^3+{q\over c\Omega(\bar{\phi},\phi)^3}\,\phi^3:\,c\in \hat{k}^*
\}
$$
and hence
$$
{\mu}^{-1}(X)=\{
c\bar{\phi}^3+{q\over c\Omega(\bar{\phi},\phi)^3}\,\phi^3:\,c\in \hat{k}^*,\,\bar{c}=
{q\over c\Omega(\bar{\phi},\phi)^3}
\}.
$$
This together with (i) implies (ii).
\end{demo}

\vskip0,5cm
\centerline{\bf Orbit parameters for ${\cal O}({\hat{k}})$}

It is clear that ${\cal O}({\hat{k}})$ is stable under the action of $\textrm{Sl}(2,k)$ and  in this section we will give a parametrisation of the space of orbits. 

Let $P\in{\cal O}({\hat{k}})$. Then, since $Q_n(P)\in \hat{k}^*{^2}$, the $\textrm{Sl}(2,\hat{k})$ orbit of $P$
regarded  as a binary cubic over $\hat{k}$ is entirely determined 
by ${I_{\widehat{\cal O}_{[1]}}}(P)$ where  
$$
{I_{\widehat{\cal O}_{[1]}}}:\widehat{\cal O}_{[1]}\rightarrow\hat{k}^*\times_{Z_2}{\hat{k}^*/\,\hat{k}^*{^3}}
$$
 is the $\textrm{Sl}(2,\hat{k})$-invariant function defined in Theorem \ref{spofinv}.  Recall that to calculate ${I_{\widehat{\cal O}_{[1]}}}(P)$,  we choose  $\lambda\in\hat{k}^*$ and $\alpha\in{\hat{k}^2}{^*}$ such that 
$$
P=\lambda\alpha^3+\bar{\lambda}\bar{\alpha}^3
$$
and then by definition,
\begin{equation}\label{hatinvdef}
{I_{\widehat{\cal O}_{[1]}}}(P)=[\hat{\omega}(\lambda\alpha^3,\overline{\lambda}\overline{\alpha}),[{\lambda}\bar{\lambda}^{-1}]].
\end{equation}
The square roots $\pm{\hat{\omega}}(\lambda\alpha^3,\overline{\lambda}\overline{\alpha}^3)$ of $Q_n(P)$ are  pure imaginary since
$$
\overline{\hat{\omega}(\lambda\alpha^3,\overline{\lambda_1}\overline{\alpha}^3)}=
\hat{\omega}(\overline{\lambda}\overline{\alpha}^3,\lambda\alpha^3)=
-\hat{\omega}(\lambda\alpha^3,\overline{\lambda_1}\overline{\alpha}^3),
$$
and  the class $[{\lambda}\bar{\lambda}^{-1}]$ of ${\lambda}\bar{\lambda}^{-1}$ in the group
${\hat{k}^*/\,\hat{k}^*{^3}}$ satisfies
$$
[{\lambda}\bar{\lambda}^{-1}]\,\overline{[{\lambda}\bar{\lambda}^{-1}]}=1.
$$
It follows
 that
$$
{I_{\widehat{\cal O}_{[1]}}}(P)\in {\rm Im}\,{\hat{k}}^*\times_{Z_2} U({\hat{k}^*/\,\hat{k}^*{^3}})
$$
where 
$$
U({\hat{k}^*/\,\hat{k}^*{^3}})=\{\alpha\in{\hat{k}^*/\,\hat{k}^*{^3}}\text{ s.t. }\alpha\bar{\alpha}=1\}.
$$
is the `unitary' group of $\hat{k}^*/\,\hat{k}^*{^3}$. Note that the $Z_2$ action on ${\rm Im}\,{\hat{k}}^*\times U({\hat{k}^*/\,\hat{k}^*{^3}})$ is precisely the natural action of ${\rm Gal}(\hat{k}/k)$.

\begin{theo}
Let  ${{I_{\widehat{\cal O}_{[1]}}}}:{\cal O}({\hat{k}})
\rightarrow
{\rm Im}\,{\hat{k}}^*\times_{ Z_2}U({\hat{k}^*/\,\hat{k}^*{^3}})$ be
defined by  \eqref{hatinvdef} above.
\vskip0,1cm
\noindent(i) Let $P,P'\in {\cal O}({\hat{k}})$. Then 
$$ 
\textrm{Sl}(2,k)\cdot P'= \textrm{Sl}(2,k)\cdot P \iff
{{I_{\widehat{\cal O}_{[1]}}}}(P')={{I_{\widehat{\cal O}_{[1]}}}}(P).
$$

\noindent(ii) The map ${{I_{\widehat{\cal O}_{[1]}}}}$ induces a bijection  
$$
{{\cal O}({\hat{k}})}/\textrm{Sl}(2,k)\longleftrightarrow
{\rm Im}\,{\hat{k}}^*\times_{Z_2}U({\hat{k}^*/\,\hat{k}^*{^3}}).
$$
\noindent(iii) The isotropy group of $P\in{\cal O}({\hat{k}})$ is isomorphic to
$$
\{
\begin{pmatrix}
\lambda&0\\
0&\bar{\lambda}
\end{pmatrix}
\in\textrm{Sl}(2,\hat{k}):\,
\lambda^3=1,\,\lambda\bar{\lambda}=1 \}.
$$
\end{theo}
\begin{demo} (i): The function
${{I_{\widehat{\cal O}_{[1]}}}}:{\cal O}({\hat{k}})
\rightarrow
{\rm Im}\,{\hat{k}}^*\times_{ Z_2}U({\hat{k}^*/\,\hat{k}^*{^3}})$ is $\textrm{Sl}(2,k)$-invariant
since it is by definition the restriction of an $\textrm{Sl}(2,\hat{k})$-invariant function on  a larger space.

To prove ${{I_{\widehat{\cal O}_{[1]}}}}$ separates orbits, suppose
${{I_{\widehat{\cal O}_{[1]}}}}(P')={{I_{\widehat{\cal O}_{[1]}}}}(P)$. Writing
$
P=\lambda\alpha^3+\bar{\lambda}\overline{\alpha}^3
$
and
$
P'=\lambda'\alpha'^3+\bar{\lambda'}\overline{\alpha'}^3,
$
 there exists $\sigma\in Z_2$ such that
\begin{equation}\label{convcond}
( \, \hat{\omega}(\lambda'\alpha'^3,\,\bar{\lambda'}\overline{\alpha'}^3),\,[\lambda'\bar{\lambda'}^{-1}]\,)
=\sigma\cdot
( \, \hat{\omega}(\lambda\alpha^3,\,\bar{\lambda}\overline{\alpha}^3),\,[{\lambda}\bar{\lambda}^{-1}]\,)
\end{equation}
and, permuting cube terms if necessary, we can suppose  without loss of generality that $\sigma$ is the identity.
Then,
equation \eqref{convcond} implies
\begin{equation}\label{convcond1}
\hat{\omega}(\lambda'\alpha'^3,\,\bar{\lambda'}\overline{\alpha'}^3)=
\hat{\omega}(\lambda\alpha^3,\,\bar{\lambda}\overline{\alpha}^3),
\quad
[\lambda'\bar{\lambda'}^{-1}]=[{\lambda}\bar{\lambda}^{-1}]
\end{equation}
or equivalently,
$$
\lambda'\bar{\lambda'}\hat{\omega}(\alpha'^3,\,\overline{\alpha'}^3)=
\lambda\bar{\lambda}\hat{\omega}(\alpha^3,\,\overline{\alpha}^3),
\quad
[\lambda'\bar{\lambda'}^{-1}]=[{\lambda}\bar{\lambda}^{-1}]
$$
which  by \eqref{sympcube} is equivalent to
\begin{equation}\label{cuberel}
\lambda'\bar{\lambda'}\,\widehat{\Omega}(\alpha',\,\overline{\alpha'})^3=
\lambda\bar{\lambda}\,\widehat{\Omega}(\alpha,\,\overline{\alpha})^3,
\quad
[\lambda'\bar{\lambda'}^{-1}]=[{\lambda}\bar{\lambda}^{-1}].
\end{equation}
Taking classes in ${\hat{k}^*/\,\hat{k}^*{^3}}$ we get
$$
[\lambda'\bar{\lambda'}]=
[\lambda\bar{\lambda}],
\quad
[\lambda'\bar{\lambda'}^{-1}]=[{\lambda}\bar{\lambda}^{-1}]
$$
and multiplying the two  equations gives
$$
[{\lambda'}^2]=[{\lambda}^2].
$$
From this it follows that $[{\lambda'}]=[{\lambda}]$ since the cube of any element in  ${\hat{k}^*/\,\hat{k}^*{^3}}$ is the identity.

Let now $\xi\in\hat{k}^*$ be such that 
$$
\lambda'=\xi^3\lambda.
$$
Substituting in  the first equation of \eqref{cuberel} we get
$$
\widehat{\Omega}(\xi\alpha',\,\overline{\xi\alpha'})^3=
\widehat{\Omega}(\alpha,\,\overline{\alpha})^3
$$
which means
$$
\widehat{\Omega}(\xi\alpha',\,\overline{\xi\alpha'})=
j\,\widehat{\Omega}(\alpha,\,\overline{\alpha})
$$
for some $j\in\hat{k}$ such that $j^3=1$. The conjugate of this equation is
$$
-\widehat{\Omega}(\xi\alpha',\,\overline{\xi\alpha'})=
-\bar{j}\,\widehat{\Omega}(\alpha,\,\overline{\alpha})
$$
and hence $\bar{j}=j$.

Define $g\in\textrm{Gl}(2,\hat{k})$ by
$$
g\cdot \alpha=j\xi\alpha',\quad g\cdot \overline{\alpha}=j\bar{\xi}\,\overline{\alpha'}.
$$
Then $g$ commutes with conjugation by definition, and preserves $\widehat{\Omega}$ since
$$
\widehat{\Omega}(g\cdot \alpha,g\cdot \overline{\alpha})=j^2\widehat{\Omega}(\xi\alpha',\,\overline{\xi\alpha'})=
j^3\widehat{\Omega}(\alpha,\overline{\alpha})=
\widehat{\Omega}(\alpha,\overline{\alpha}).
$$
Hence $g\in \textrm{Sl}(2,k)$. Furthermore,
$$
g\cdot P=\lambda(g\cdot\alpha)^3+\bar{\lambda}(\overline{g\cdot \alpha})^3
=\lambda(j\xi\alpha')^3+\bar{\lambda}(j\bar{\xi}\,\overline{\alpha'})^3
=\lambda'\alpha'^3+\bar{\lambda'}\overline{\alpha'}^3=P'
$$
which shows  that $P$ and $P'$ are in the same $\textrm{Sl}(2,k)$-orbit. This proves (i).

To prove (ii), we only have  to show that ${{I_{\widehat{\cal O}_{[1]}}}}$ is surjective since by (i), the function
${{I_{\widehat{\cal O}_{[1]}}}}$ induces an injection ${{\cal O}({\hat{k}})}/\textrm{Sl}(2,k)\hookrightarrow
{\rm Im}\,{\hat{k}}^*\times_{ Z_2}U({\hat{k}^*/\,\hat{k}^*{^3}})$. 

Let $(q,s)\in
{\rm Im}\,{\hat{k}}^*\times U({\hat{k}^*/\,\hat{k}^*{^3}})$. 
First, pick $\lambda\in{\hat{k}}^*$ such that 
 \begin{equation}\label{surj0}
[\lambda]=\bar{s}.
 \end{equation}
 Since $[\lambda\bar{\lambda}]=s\bar{s}=1$, we know $\lambda\bar{\lambda}$ is a cube in $\hat{k}{^*}$ but in fact, since $\hat{k}{^*}/k$ is a quadratic extension and $\lambda\bar{\lambda}\in k$, this implies that there exists $r\in {k}{^*}$ such that 
 \begin{equation}\label{surj1}
 \lambda\bar{\lambda}=r^3.
 \end{equation}
Now let
$$
\alpha=-{q\over2r}\hat{x}+\hat{y}
$$ 
(where $\hat{x},\hat{y}\in \hat{k}^2{^*}$ are the base extensions of $x,y\in {k^2}^*$) and let
\begin{equation}\label{invimqs}
P={\lambda\over q}\alpha^3-{\bar{\lambda}\over q}\,\overline{ \alpha}^3.
\end{equation}
This is a binary cubic of the form $T+\bar{T}$ where $T\in\hat{T}$.  We are now going to show that
$P\in {\cal O}({\hat{k}})$ and that ${{I_{\widehat{\cal O}_{[1]}}}}(P)=[q,s]$.

Note first that
$$
\widehat{\Omega}(\alpha,\bar{\alpha})=
\widehat{\Omega}(-{q\over2r}\hat{x}+\hat{y},-\overline{({q\over2r})}\hat{x}+\hat{y}=-{q\over2r}+\overline{({q\over2r})}=-{q\over r},
$$
so $\widehat{\Omega}(\alpha,\bar{\alpha})\not=0$ which means $\alpha$ and $\overline{\alpha}$ are not proportional. Hence  $\alpha^3$ and $\overline{\alpha}^3$  are coprime and $P\in{ \cal O}({\hat{k}})$.

Next, we have
 \begin{equation}\label{surj2}
\hat{\omega}(\alpha^3,\bar{\alpha}^3)=\widehat{\Omega}(\alpha,\bar{\alpha})^3
=-{q^3\over r^3}
\end{equation}
 and 
 \begin{equation}\label{surj3}
\hat{\omega}({\lambda\over q}\alpha^3,-{\bar{\lambda}\over q}\overline{ \alpha}^3)=
-({1\over q})^2\lambda\bar{\lambda}\hat{\omega}(\alpha^3,\bar{\alpha}^3)=
q
\end{equation}
using equations \eqref{surj1} and \eqref{surj2}.
Finally, it follows from  \eqref{surj0} that
 \begin{equation}\label{surj4}
\left[{\lambda\over q}(\overline{\lambda\over q})^{-1}\right]=
\left[\lambda\bar{\lambda}^{-1}\right]=
\bar{s}s^{-1}=
s^{-1}s^{-1}=s^{-2}=s.
\end{equation}
Hence, putting together equations \eqref{hatinvdef}, \eqref{invimqs}, \eqref{surj3} and \eqref{surj4}, we get
$$
{{I_{\widehat{\cal O}_{[1]}}}}(P)=[q,s]
$$
and this proves  that
${{I_{\widehat{\cal O}_{[1]}}}}:{{\cal O}({\hat{k}})}\rightarrow
{\rm Im}\,{\hat{k}}^*\times_{ Z_2}U({\hat{k}^*/\,\hat{k}^*{^3}})$
is surjective.

Part (iii) follows from Theorem \ref{spofinv} (iii).

\end{demo}
\begin{cor} Let $P,P'\in {\cal O}({\hat{k}})$. Then 
$$ 
\textrm{Sl}(2,k)\cdot P'= \textrm{Sl}(2,k)\cdot P \iff
\textrm{Sl}(2,\hat{k})\cdot P'= \textrm{Sl}(2,\hat{k})\cdot P .
$$
\end{cor}
\begin{demo} Both properties are equivalent to ${{I_{\widehat{\cal O}_{[1]}}}}(P)={{I_{\widehat{\cal O}_{[1]}}}}(P')$ by the above theorem and Theorem \ref{spofinv}.
\end{demo}
\vskip 0,1cm
\centerline{\bf Properties of orbit space}
\vskip 0,1cm
The parameter space 
$$
{\rm Im}\,{\hat{k}}^*\times_{ Z_2}U({\hat{k}^*/\,\hat{k}^*{^3}})
$$
 for
$\textrm{Sl}(2,k)$ orbits in ${\cal O}({\hat{k}})$ is very analogous to the parameter space 
$$
{\hat{k}}^*\times_{ Z_2}{{k}^*/\,{k}^*{^3}}
$$
 for
$\textrm{Sl}(2,k)$ orbits in ${\cal O}_{[1]}$ that we gave in Theorem \ref{spofinv}. Its main features  can best be summarized in the diagram
\vskip0,1cm
\begin{equation}\label{unitarydiag}
 \xymatrix
{
&{\rm Im}\,{\hat{k}}^*\times_{ Z_2}U({\hat{k}^*/\,\hat{k}^*{^3}})\ar[dl]^{{\widehat{sq}}}\ar[dr]^{{\hat{t}}}&\\
({\rm Im}\,{\hat{k}}^*)^2\ar@/^/[ur]^{\hat{e}}&&U({\hat{k}^*/\,\hat{k}^*{^3}})/Z_2.\\
}
 \end{equation}
 \vskip0,1cm
The map
\begin{equation}\label{unitarysq}
{\widehat{sq}}:{\rm Im}\,{\hat{k}}^*\times_{ Z_2}U({\hat{k}^*/\,\hat{k}^*{^3}})
\rightarrow({\rm Im}\,{\hat{k}}^*)^2
\end{equation}
given by
$$
{\widehat{sq}}([q,\alpha])=q^2
$$
is the fibration associated to the principal $Z_2$-fibration
$$
{\rm Im}\,{\hat{k}^*}\rightarrow({\rm Im}\,{\hat{k}}^*)^2
$$
and the action of $Z_2$ on $U({\hat{k}^*/\,\hat{k}^*{^3}})$ by conjugation. Since $Z_2$ acts by automorphisms, the fibre ${\widehat{sq}}^{-1}(q^2)$ over any point  $q^2\in ({\rm Im}\,{\hat{k}}^*)^2$ has
a natural group structure
\begin{equation}\label{gpstrunitary}
[q,u_1]\times[q,u_2]=[q,u_1u_2]
\end{equation}
independent of the choice of square root $q$ of $q^2$. Taking the identity at each point, we get a
 canonical section 
$\hat{e}:({\rm Im}\,{\hat{k}}^*)^2\rightarrow {\rm Im}\,{\hat{k}}^*\times_{ Z_2}U({\hat{k}^*/\,\hat{k}^*{^3}})$ of 
\eqref{unitarysq} given by
\begin{equation}\label{iddiscfix}
\hat{e}(q^2)=[q,1]
\end{equation}
but, although each fibre is a group isomorphic to $U({\hat{k}^*/\,\hat{k}^*{^3}})$, the fibration
\eqref{unitarysq} is not in general isomorphic to the product 
$$
({\rm Im}\,{\hat{k}}^*)^2\times U({\hat{k}^*/\,\hat{k}^*{^3}})
\rightarrow({\rm Im}\,{\hat{k}}^*)^2.
$$

To translate the above features of orbit space into  more concrete statements about binary cubics over $k$, note that the map
$\widehat{sq}$ is essentially the quartic $Q_n$ since for all $P\in{\cal O}({\hat{k}})$,
$$
\widehat{sq}({{I_{\widehat{\cal O}_{[1]}}}}(P))=Q_n(P).
$$
\begin{theo}\label{fixeddisc}
Let $M\in ({\rm Im}\,{\hat{k}}^*)^2$, let 
$$
{\cal O}_M=\{P\in S^3({k^2}^*):\,Q_n(P)=M\}
$$
and let ${\cal O}_M/\textrm{Sl}(2,k)$ be the  space of ${Sl}(2,k)$-orbits in ${\cal O}_M$. 

\medskip
\noindent(i) The map ${I_{\widehat{\cal O}_{[1]}}}:{\cal O}({\hat{k}})\rightarrow{\rm Im}\,{\hat{k}}^*\times_{ Z_2}U({\hat{k}^*/\,\hat{k}^*{^3}})$ induces a bijection
$$
{\cal O}_M/\textrm{Sl}(2,k)\longleftrightarrow
{\widehat{sq}}^{-1}(M)
$$
and, by pullback of\eqref{gpstrunitary}, a group structure on ${\cal O}_M/\textrm{Sl}(2,k)$.

\medskip
\noindent(ii) As groups, ${\cal O}_M/\textrm{Sl}(2,k)\cong U({\hat{k}^*/\,\hat{k}^*{^3}})$.

\medskip
\noindent(iii) The identity element of ${\cal O}_M/\textrm{Sl}(2,k)$ is characterised by:
$$
\textrm{Sl}(2,k)\cdot P=1\Leftrightarrow P\text{ is reducible over }k.
$$
\end{theo} 
\begin{demo} Parts (i) and (ii) follow from the discussion above. To prove (iii), first note that
$P$ is reducible over $k$ iff $P$ is reducible over $\hat{k}$ since $P$ is  cubic and $\hat{k}/k$ is
a quadratic extension. By Theorem \ref{fixeddiscnonunitary}(iii), $P$ is reducible over $\hat{k}$ iff
${I_{\widehat{\cal O}_{[1]}}}(P)=[q,1]$  where $q\in\hat{k}$ is a square root of $M$, and 
by equation \eqref{iddiscfix}, this is the identity element of ${\cal O}_M/\textrm{Sl}(2,k)$.
\end{demo}
\begin{cor} \label{uniqredorb2}
Let $P,P'\in S^3({k^2}^*)$ be reducible binary cubics such that $Q_n(P)=Q_n(P')$
is nonzero. Then there exists $g\in\textrm{Sl}(2,k)$ such that $P'=g\cdot P$.
\end{cor}
\begin{demo} Suppose $Q_n(P)=Q_n(P')=M$. If $M\in{k^*}^2$, the result follows from Theorem \ref{fixeddiscnonunitary}(iii). If $M\in{k^*}$  is not a square, one can find a quadratic extension $\hat{k}$ of $k$ such that  $M\in ({\rm Im}\,{\hat{k}}^*)^2$. The result then follows from Theorem \ref{fixeddisc} (iii).
\end{demo}

To finish this section we briefly discuss the map
${\hat{t}}:{\rm Im}\,{\hat{k}}^*\times_{ Z_2}U({\hat{k}^*/\,\hat{k}^*{^3}})\rightarrow U({\hat{k}^*/\,\hat{k}^*{^3}})/Z_2$ in diagram \eqref{unitarydiag}
given by
$$
{\hat{t}}([q,\alpha])=[\alpha].
$$
This a fibration with fibre ${\rm Im}\,{\hat{k}}^*$ outside the identity coset $[1]$ but  
 $$
 {\hat{t}}^{-1}([1])=\hat{e}({k^*}^2)
 $$
  is a `singular fibre'.
   There is a $k^*$-action:
\begin{equation}
\lambda\cdot [q,\alpha]=[\lambda q,\alpha]
\end{equation}
which maps fibres of ${\widehat{sq}}$ to fibres of ${\widehat{sq}}$:
$$
{\widehat{sq}}([q',\alpha'])={\widehat{sq}}([q,\alpha])\Rightarrow {\widehat{sq}}(\lambda\cdot[q',\alpha'])={\widehat{sq}}(\lambda\cdot[q,\alpha]),
$$
and whose orbits are exactly the fibres of ${\hat{t}}$:
$$
{\hat{t}}([q',\alpha'])=t([q,\alpha])\Leftrightarrow \exists\lambda\in k^*\text{ s.t. }
[q',\alpha']=\lambda\cdot [q,\alpha].
$$
Isotropy for this action is given by: 
$
Isot_{k^*}([q,\alpha])=\left\{
\begin{array}{lr}
1&\text{if }\alpha\not=1\\
\{\pm 1\}&\text{if }\alpha=1.
\end{array}
\right.
$

It would be interesting to interpret these features of the orbit space in terms of the original binary cubics.

\section{ Parameter spaces for $\textrm{Gl}(2,k)$-orbits }

  We have seen that the $\textrm{Sl}(2,k)$-orbits in 
 $$
{\cal O}_{[1]}=\{P\in S^3({k^2}^*):\,Q_n(P)\in{k^*}^2\}
 $$ are parametrised by
 $$
 k^*\times_{Z_2}k^*/{k^*}^3
 $$
 and that  if $\hat{k}$ is a quadratic extension of $k$, the $\textrm{Sl}(2,k)$-orbits in 
 $$
 {\cal O}(\hat{k})=\{P\in S^3({k^2}^*):\,Q_n(P)\in{({\rm Im}\,\hat{k}^*)}^2\}
 $$ are parametrised by
 $$
{\rm Im}\,\hat{k}^*\times_{Z_2}U({\hat{k}}^*/  {\hat{k}}^{* 3 }).
 $$
 The group $GL(2,k)$ also acts on binary cubics and since
 $$
 Q_n(g\cdot P)=({\rm det}\,g)^{-6}Q_n(P)\quad\forall g\in \textrm{Gl}(2,k),\forall P\in S^3({k^2}^*),
 $$
 the spaces ${\cal O}_{[1]}$ and $ {\cal O}(\hat{k})$ are stable under $\textrm{Gl}(2,k)$.

 In general, if $\textrm{Gl}(2,k)$ acts on a space $X$ there is a map
$$
X/\textrm{Sl}(2,k)\rightarrow X/\textrm{Gl}(2,k)
$$
from the set of $\textrm{Sl}(2,k)$-orbits onto the set of $\textrm{G}(2,k)$-orbits. The fibres of
this map are the orbits of the $k^*$-action on $X/\textrm{Sl}(2,k)$ given by
\begin{equation}\label{SLorbtoGLorb}
\lambda* [x]=[\Lambda\cdot x]
\end{equation}
where $\Lambda$ is any element of $\textrm{Gl}(2,k)$ such that ${\rm det}\,\Lambda=\lambda$.
Thus to get parameter spaces for ${\cal O}_{[1]}/\textrm{Gl}(2,k)$ and
$ {\cal O}(\hat{k})/\textrm{Gl}(2,k)$ we need just to calculate the $k^*$-actions on
 $k^*\times_{Z_2}k^*/{k^*}^3$ and ${\rm Im}\,\hat{k}^*\times_{Z_2}U({\hat{k}}^*/  {\hat{k}}^{* 3 })$
 corresponding to \eqref{SLorbtoGLorb}.
\begin{lem} \label{lawunderGL}
(i) Let $k^*$ act on $k^*\times_{Z_2}k^*/{k^*}^3$ by
$$
\lambda\cdot[\xi,\alpha]=[\lambda\xi,\alpha]
$$
and let $I_{{\cal O}_{[1]}}:{\cal O}_{[1]}\rightarrow k^*\times_{Z_2}k^*/{k^*}^3$ be defined
by \eqref{diaginv}. Then
$$
I_{{\cal O}_{[1]}}(g\cdot P)=({\rm det}\,g)^{-3}\cdot I_{{\cal O}_{[1]}}(P)
\qquad\forall P\in{\cal O}_{[1]},\,\forall g\in\textrm{Gl}(2,k).
$$

\noindent(ii) Let $k^*$ act on ${\rm Im}\,\hat{k}^*\times_{Z_2}U({\hat{k}}^*/  {\hat{k}}^{* 3 })$ by
$$
\lambda\cdot[\xi,\alpha]=[\lambda\xi,\alpha]
$$
and let $I_{\widehat{\cal O}_{[1]}}: {\cal O}(\hat{k})\rightarrow {\rm Im}\,\hat{k}^*\times_{Z_2}U({\hat{k}}^*/  {\hat{k}}^{* 3 })$ be defined
by \eqref{hatinvdef}. Then
$$
I_{\widehat{\cal O}_{[1]}}(g\cdot P)=({\rm det}\,g)^{-3}\cdot I_{{\cal O}_{[1]}}(P)
\qquad\forall P\in{\cal O}(\hat{k}),\,\forall g\in\textrm{Gl}(2,k).
$$
\end{lem}
\begin{demo} To prove (i), since for any $P\in {\cal O}_{[1]}$ there exists $h\in\textrm{Gl}(2,k)$ and $a,b\in k^*$ such that
$$
h\cdot P=ax^3+by^3,
$$
 it is sufficient to prove that
 $$
 I_{{\cal O}_{[1]}}(g\cdot (ax^3+by^3))=({\rm det}\,g)^{-3}\cdot I_{{\cal O}_{[1]}}(ax^3+by^3)
\qquad\forall a,b\in k^*,\,\forall g\in\textrm{Gl}(2,k).
 $$
 Consider
 $$
 g'=
 \begin{pmatrix}
 {\rm det}\,g&0\\
 0&1
 \end{pmatrix}.
 $$
 Then $ {\rm det}\,g'= {\rm det}\,g$, $g'\cdot x={1\over {\rm det}\,g}x$ and $g'\cdot y=y$. Hence
 $g'g^{-1}\in\textrm{Sl}(2,k)$,
 $$
 I_{{\cal O}_{[1]}}(g\cdot (ax^3+by^3))=
   I_{{\cal O}_{[1]}}(g'\cdot (ax^3+by^3))
 $$
 and
 $$
  I_{{\cal O}_{[1]}}(g'\cdot (ax^3+by^3))= I_{{\cal O}_{[1]}}({a\over({\rm det}\,g)^3}x^3+by^3))
  =[{ab\over({\rm det}\,g)^3},[ab^{-1}]].
 $$
 The result follows since $ I_{{\cal O}_{[1]}}(ax^3+by^3)=[ab,[ab^{-1}]]$.
 
 Part (ii) follows  from (i) applied to $\hat{k}$.
\end{demo}

\begin{cor} 
(i) If $P\in{\cal O}_{[1]}$ and $\lambda\in k^*$ then
$$
I_{{\cal O}_{[1]}}(\lambda*[P])={1\over\lambda^3}\cdot I_{{\cal O}_{[1]}}([P]).
$$

\noindent(ii) If $P\in{\cal O}(\hat{k})$ and $\lambda\in k^*$ then
$$
I_{\widehat{\cal O}_{[1]}}(\lambda*[P])={1\over\lambda^3}\cdot I_{\widehat{\cal O}_{[1]}}([P]).
$$

\end{cor}
\begin{demo} Immediate from the lemma.
\end{demo}

From this we get the $k^*$-actions on the parameter spaces 
 $k^*\times_{Z_2}k^*/{k^*}^3$ and ${\rm Im}\,\hat{k}^*\times_{Z_2}U({\hat{k}}^*/  {\hat{k}}^{* 3 })$
 corresponding to \eqref{SLorbtoGLorb} : $\lambda\in k^*$ acts by multiplication by $\lambda^{-3}$ on the first factor.  
 
 Hence, by the discussion above, the maps $I_{{\cal O}_{[1]}}:{\cal O}_{[1]}\rightarrow k^*\times_{Z_2}k^*/{k^*}^3$  and
$I_{\widehat{\cal O}_{[1]}}: {\cal O}(\hat{k})\rightarrow {\rm Im}\,\hat{k}^*\times_{Z_2}U({\hat{k}}^*/  {\hat{k}}^{* 3 })$ induce  bijections
$$
\begin{array}{rll}
{\cal O}_{[1]}/\textrm{Gl}(2,k)&\longleftrightarrow  (k^*\times_{Z_2}k^*/{k^*}^3)/k^{*3}&=
k^*/k^{*3}\times\, (k^*/k^{*3})/{Z_2},\\
{\cal O}(\hat{k})/\textrm{Gl}(2,k)&\longleftrightarrow  ( {\rm Im}\,\hat{k}^*\times_{Z_2}U({\hat{k}}^*/  {\hat{k}}^{* 3 }))/k^{*3}&=({\rm Im}\,\hat{k}^*)/k^{*3}\times \,U({\hat{k}}^*/  {\hat{k}}^{* 3 })/Z_2.
\end{array}
$$
To summarize, we have proved the
  \begin{theo} \label{parGLorb}
  (a) Define $\pi: k^*\times_{Z_2}k^*/{k^*}^3\rightarrow k^*/k^{*3}\times\, (k^*/k^{*3})/{Z_2}$ by $\pi([\xi,\alpha])=([\xi],[\alpha])$ and $J_{{\cal O}_{[1]}}:{\cal O}_{[1]}\rightarrow k^*/k^{*3}\times\, (k^*/k^{*3})/{Z_2}$ by $J_{{\cal O}_{[1]}}=\pi\circ  I_{{\cal O}_{[1]}}$. 
  
  (i) Let $P,P'\in{\cal O}_{[1]}$ Then
 $$
 \textrm{Gl}(2,k)\cdot P= \textrm{Gl}(2,k)\cdot P'\quad\Leftrightarrow\quad J_{{\cal O}_{[1]}}(P)=J_{{\cal O}_{[1]}}(P').
 $$

(ii) The map $J_{{\cal O}_{[1]}}$ induces a bijection
 $$
 {\cal O}_{[1]}/\textrm{Gl}(2,k)\longleftrightarrow k^*/k^{*3}\times\, (k^*/k^{*3})/{Z_2}.
 $$
 
\noindent (b) Define $\hat{\pi}: {\rm Im}\,\hat{k}^*\times_{Z_2}U({\hat{k}}^*/  {\hat{k}}^{* 3 })\rightarrow ({\rm Im}\,\hat{k}^*)/k^{*3}\times \,U({\hat{k}}^*/  {\hat{k}}^{* 3 })/Z_2$ by $\hat{\pi}([\xi,\alpha])=([\xi],[\alpha])$ and $J_{\widehat{\cal O}_{[1]}}:{\cal O}(\hat{k})\rightarrow k^*/k^{*3}\times\, (k^*/k^{*3})/{Z_2}$ by $J_{\widehat{\cal O}_{[1]}}=\hat{\pi}\circ  I_{\widehat{\cal O}_{[1]}}$.

(i)  Let
 $P,P'\in{\cal O}(\hat{k})$. Then
 $$
 \textrm{Gl}(2,k)\cdot P= \textrm{Gl}(2,k)\cdot P'\quad\Leftrightarrow\quad J_{\widehat{\cal O}_{[1]}}(P)=J_{\widehat{\cal O}_{[1]}}(P').
 $$
 
 (ii) The map $J_{\widehat{\cal O}_{[1]}}$ induces a bijection
 $$
 {\cal O}(\hat{k})/\textrm{Gl}(2,k)\longleftrightarrow ({\rm Im}\,\hat{k}^*)/k^{*3}\times \,U({\hat{k}}^*/  {\hat{k}}^{* 3 })/Z_2.
 $$
 \end{theo}
 \vskip0,1cm
\centerline{\bf{Orbits spaces when $k$ is a finite field of characteristic not $2$ or $3$}}
 \vskip0,1cm
Let $k$ be a finite field with $q$ elements, not of characteristic  $2$ or $3$. The following facts 
are well-known:
\begin{itemize}
\item
$k^*/k^{*2}\cong Z_2$ so up to isomorphsm, there is only one quadratic extension of $k$ and
$k^{*2}$ has ${1\over2}(q-1)$ elements;
\item
if $q=1$ mod $3$, $k^*/k^{*3}\cong \bb Z/3\bb Z$;
\item
if $q=2$ mod $3$, $k^*=k^{*3}$;
\item
if $q=1$ mod $3$ and $\hat{k}/k$ is a quadratic extension, $U({\hat k}^*/{\hat k}^{*3})\cong 1$;
\item
if $q=2$ mod $3$ and $\hat{k}/k$ is a quadratic extension, $U({\hat k}^*/{\hat k}^{*3})\cong \bb Z/3\bb Z$.
\end{itemize}
These facts together with Theorem \eqref{spofinv} and Theorem \ref{parGLorb} immediately give the
\begin{prop}\label{orbcountfinfield}
Let $k$ be a finite field with $q$ elements, not of characteristic  $2$ or $3$
and let $\hat k$ be a quadratic extension. Set
$$
\begin{array}{l}
{\cal O}_{[1]}=\{P\in S^3({k^2}^*):\,Q_n(P)\in{k^*}^2\},\\
 {\cal O}(\hat{k})=\{P\in S^3({k^2}^*):\,Q_n(P)\in{({\rm Im}\,\hat{k}^*)}^2\}.
\end{array}
$$

(a) If $q=1$ mod 3, ${\cal O}_{[1]}$ is  the union of ${3\over 2}(q-1)$ $\textrm{Sl}(2,k)$-orbits and
$ {\cal O}(\hat{k})$ is  the union of ${1\over 2}(q-1)$ $\textrm{Sl}(2,k)$-orbits.

\medskip
(b)  If $q=1$ mod 3, ${\cal O}_{[1]}$ is  the union of $6$ $\textrm{Gl}(2,k)$-orbits and
$ {\cal O}(\hat{k})$ is  the union of $3$ $\textrm{Gl}(2,k)$-orbits.

\medskip
(c) If $q=2$ mod 3, ${\cal O}_{[1]}$ is  the union of ${1\over 2}(q-1)$ $\textrm{Sl}(2,k)$-orbits and
$ {\cal O}(\hat{k})$ is  the union of ${3\over 2}(q-1)$ $\textrm{Sl}(2,k)$-orbits.

\medskip
(d)  If $q=2$ mod 3, ${\cal O}_{[1]}$ is  a $\textrm{Gl}(2,k)$-orbit and
$ {\cal O}(\hat{k})$ is  the union of $2$ $\textrm{Gl}(2,k)$-orbits.
\end{prop}
\begin{demo} As examples, let us count the number of $\textrm{Sl}(2,k)$-orbits in ${\cal O}_{[1]}$ when $q=1$ mod 3 and the number of $\textrm{Gl}(2,k)$-orbits in ${\cal O}(\hat{k})$ when $q=2$ mod 3.

In the first case, by Theorem \eqref{spofinv}, the parameter space is  $k^*\times_{Z_2}k^*/{k^*}^3$ which, being a fibre bundle over $k^{*2}$ with fibre $k^*/k^{*3}$, has ${1\over 2}(q-1)\times 3={3\over 2}(q-1)$ elements. 

In the second case, by Theorem \ref{parGLorb} , the parameter space is  
$({\rm Im}\,\hat{k}^*)/k^{*3}\times \,U({\hat{k}}^*/  {\hat{k}}^{* 3 })/Z_2$
and this has $1\times 2=2$ elements since $Z_2$ acts on $U({\hat{k}}^*/  {\hat{k}}^{* 3 })$ by
inversion.
\end{demo}
According to \cite{HM} (Proposition 5.6) at least part of the following corollary can be found in Dickson \cite{DiMad}.
\begin{cor}
Let $k$ be a finite field with $q$ elements, not of characteristic  $2$ or $3$. The number of
$\textrm{Sl}(2,k)$-orbits of binary cubics with nonzero discriminant is $2(q-1)$. The number of
$\textrm{Gl}(2,k)$-orbits of binary cubics with nonzero discriminant is $9$ if $q=1$ mod 3 and 
$3$ if $q=2$ mod 3.
\end{cor}
\begin{demo}
A binary cubic of nonzero discriminant is either in ${\cal O}_{[1]}$ or in
${\cal O}(\hat{k})$ since up to isomorphism, $k$ has only one quadratic extension. Hence, the total
number of $\textrm{Sl}(2,k)$-orbits with nonzero discriminant is the number of $\textrm{Sl}(2,k)$-orbits in ${\cal O}_{[1]}$ plus  the number of $\textrm{Sl}(2,k)$-orbits in ${\cal O}(\hat{k})$. The
same is true for $\textrm{Gl}(2,k)$-orbits and the result follows from Proposition \ref{orbcountfinfield}.
\end{demo}

\section{A symplectic Eisenstein identity}

The following identity is a symplectic generalisation of the classical Eisenstein identity which,
 as we will see, is obtained from it  in the special case when $Q$ is the cube of a linear form. There is an analogous identity for the symplectic module associated to any Heisenberg graded Lie algebra (\cite{Sl-St}).
\begin{theo} Let $P,Q\in S^3(k^2{^*})$. Then
\begin{align}\label{sympEisenstein}
\nonumber
\omega(\Psi(P),Q)^2
-9Q_n(P)\,\omega(P,Q)^2=&\\
-{9\over 2}\,\omega(\mu(P)^{\otimes3}\,\cdot Q,Q)&
-{9\over 2}\,Q_n(P)\,\omega(\mu(P)\cdot Q,Q)
\end{align}
where $\mu(P)^{\otimes3}$ denotes the unique endomorphism of $S^3(k^2{^*})$ satisfying
$\mu(P)^{\otimes3}\cdot(\alpha^3)=(\mu(P)\cdot\alpha)^3$ for all $\alpha\in k^2{^*}$.
\end{theo}
\begin{demo} If $\mu(P)=0$, then $\Psi(P)=0,\,Q_n(P)=0$ and all terms in the identity are zero.

If $\mu(P)$ is nilpotent nonzero, then $Q_n(P)=0$ and and there exists $g\in{\rm Sl}(2,k)$ such that $g\cdot P=x^2y$. Since the  identity  \eqref{sympEisenstein} is  ${\rm Sl}(2,k)$-invariant,  we can suppose without loss of generality that $P=x^2y$. Then, by calculation,
$$
\Psi(P)=-{2\over 9}x^3,\quad
\mu(P)=
{2\over 9}
\begin{pmatrix}
0&0\\
1&0
\end{pmatrix}
$$
and so $\mu(P)\cdot x=0$ and $\mu(P)\cdot y=-{2\over 9}x$. Let
$$
Q=px^3+3rx^2y+3sxy^2+ty^3.
$$
The LHS of \eqref{sympEisenstein} is
$$
\omega(-{2\over 9}x^3,Q)^2=({2\over 9})^2t^2.
$$
and the RHS of \eqref{sympEisenstein}  is
$$
-{9\over 2}\omega(\mu(P)^{\otimes3}\,\cdot Q,Q)=
-{9\over 2}\omega(-({2\over 9})^3tx^3,Q)=
({2\over 9})^2t^2.
$$
Thus \eqref{sympEisenstein} holds if
$\mu(P)$ is nilpotent nonzero.

To complete the proof of the proposition it remains to prove \eqref{sympEisenstein} if $Q_n(P)\not=0$.
As the identity is independent of the field we may suppose that $Q_n(P)$ is a square in $k^*$ and hence that $P\in{\cal O}_{[1]}$. Since the  identity  \eqref{sympEisenstein} is  ${\rm Sl}(2,k)$-invariant,  we can further suppose without loss of generality that 
$$
P=ax^3+dy^3.
$$
Then
$$
Q_n(P)=a^2d^2,\quad
\Psi(P)=3ad(-ax^3+dy^3),\quad
\mu(P)=
\begin{pmatrix}
ad&0\\
0&-ad
\end{pmatrix}
$$
and so $\mu(P)\cdot x=-adx$ and $\mu(P)\cdot y=ady$. Let
$$
Q=px^3+3rx^2y+3sxy^2+ty^3.
$$
The LHS of \eqref{sympEisenstein} is
\begin{align}\label{sympEisleft}
\nonumber
\omega(\Psi(P)^2,Q)&-9Q_n(P)\omega(P,Q)^2\\
\nonumber
&=9a^2d^2\big(\omega(-ax^3+dy^3,Q)^2-\omega(ax^3+dy^3,Q)^2\bigr)\\
\nonumber
&=-36a^3d^3\omega(x^3,Q)\omega(y^3,Q)\\
&=36a^3d^3pt.
\end{align}
On the other hand,  the first term of the RHS of \eqref{sympEisenstein} is
\begin{align}\label{Eisright1}
\nonumber
-{9\over 2}\omega(\mu(P)^{\otimes3}\,\cdot Q,Q)
&=
-{9\over 2}a^3d^3\omega(-px^3+3rx^2y-3sxy^2+ty^3,Q)\\
\nonumber
&=-{9\over 2}a^3d^3(-2pt-6rs)\\
&=9a^3d^3(pt+3rs)
\end{align}
and the second term of the RHS of \eqref{sympEisenstein} is
\begin{align}\label{Eisright2}
\nonumber
-{9\over 2}Q_n(P)\omega(\mu(P)\cdot Q,Q)
&=-{9\over 2}a^3d^3\omega(-3px^3-3rx^2y+3sxy^2+3ty^3,Q)\\
\nonumber
&=-{9\over 2}a^3d^3(-6pt+6rs)\\
&=27a^3d^3(pt-rs).
\end{align}
The result follows from equations \eqref{sympEisleft}, \eqref{Eisright1} and \eqref{Eisright2}.
\end{demo}

To obtain the classical Eisenstein identity from this result, recall that
one can use the symplectic form $\Omega$ on $k^2{^*}$ to define a ${\rm Sl}(2,k)$-equivariant isomorphism $\,\tilde{}:k^2\rightarrow k^2{^*}$: if $v\in k^2$, we let $\tilde{v}\in k^2{^*}$  be the unique linear form such that
$$
\phi(v)=\Omega(\phi,\tilde{v})\quad\forall\phi\in k^2{^*}.
$$
It then follows that
\begin{equation}\label{evalissymp}
P(v)=\omega(P,\tilde{v}^3)\quad\forall P\in S^3(k^2{^*}), \,\forall v\in k^2,
\end{equation}
so that  the operation of evaluating a binary cubic at a point of  $k^2$ can be  expressed in terms of the symplectic form $\omega$ on $S^3(k^2{^*})$. One can also pullback $\Omega$ to get an ${\rm Sl}(2,k)$-invariant symplectic form 
$\Omega_{k^2}$ on $k^2$:
$$
\Omega_{k^2}(v,w)=\Omega(\tilde{v},\tilde{w})\quad\forall v,w\in k^2.
$$
\begin{cor}
(Classical Eisenstein identity)
Let $P\in S^3(k^2{^*})$ and let $v\in k^2$.
\begin{equation}
\nonumber
\Psi(P)(v)^2
-9Q_n(P)\,P(v)^2=\\
-{9\over 2}\,\Omega_{ k^2}(\mu(P)\cdot v,v)^3.
\end{equation}
\end{cor}
\begin{demo} Setting $Q=\tilde{v}^3$ in \eqref{sympEisenstein} and using \eqref{evalissymp}, we get
\begin{align}
\nonumber
\Psi(P)(v)^2
-9Q_n(P)\,P(v)^2=&\\
-{9\over 2}\,\omega(\mu(P)^{\otimes3}\,\cdot \tilde{v}^3,\tilde{v}^3)&
-{9\over 2}\,Q_n(P)\,\omega(\mu(P)\cdot \tilde{v}^3,\tilde{v}^3).
\end{align}
The result follows from this since
$$
\omega(\mu(P)\cdot \tilde{v}^3,\tilde{v}^3)=
3\omega((\mu(P)\cdot \tilde{v})\tilde{v}^2,\tilde{v}^3)=
0
$$
( $(\mu(P)\cdot \tilde{v})\tilde{v}^2$ has at least a double root at $v$) and
$$
\omega(\mu(P)^{\otimes3}\,\cdot \tilde{v}^3,\tilde{v}^3)=
\Omega(\mu(P)\cdot\tilde{v},\tilde{v})^3=
\Omega_{ k^2}(\mu(P)\cdot v,v)^3.
$$
\end{demo}


\begin{thebibliography}{999}

\bibitem{Bh}{\sc M. Bhargava}, Higher composition laws. I. A new view on Gauss 
composition, and quadratic generalizations,  {\it Ann. of Math. (2) } {\bf  159}   
(2004),  no. 1, 217--250.

 \bibitem{Bour}{\sc N. Bourbaki}, {\bf   \'Elements de math\'ematique}, Fascicule XXIV, Livre II, Alg\`ebre Chapitre 9, Hermann, Paris, 1959.

 \bibitem{Ca-Schw}
{\sc M. Cahen and L. Schwachh\" ofer}, Special symplectic connections and Poisson geometry. 
{\it Lett. Math. Phys. } {\bf 69} (2004), 115--137. 


\bibitem{Dema}{\sc M. Demazure}, Automorphismes et d\'eformations des vari\'et\'es de Borel.  {\it Invent. Math. }{\bf  39 } (1977), no. 2, 179--186. 

\bibitem{DiHst}{\sc L. E. Dickson}, {\bf History of the Theory of Numbers}, vol. III. Chelsea, 1952.

\bibitem{Di} {\sc L. E. Dickson}, {\bf Algebraic Invariants}, Mathematical Monographs, 
No. 14, John Wiley and Sons, New York, 1914.

\bibitem{DiMad} {\sc L. E. Dickson}, {\bf On Invariants and The Theory of Numbers}, The Madison Colloquim, 
Dover Publications, New York, 1966.

\bibitem{Ei}{\sc G. Eisenstein}, Untersuchungen \" uber die cubischen Formen mit zwei 
Variabeln, {\it  J. Crelle } {\bf 27} (1844), 89--104 = Mathematische Werke, Band I, 
Chelsea Publ. Co., 1975, 10--25.

\bibitem{Ha}{\sc D. Haile}, On the Clifford algebra of a binary cubic form,  {\it  
Amer. J. Math. } {\bf 106} (1984), 1269--1280.  

\bibitem{HM}{\sc J. W. Hoffman and J.Morales}, Arithmetic of binary cubic forms,  {\it 
 Enseign. Math. (2)}, {\bf 46}  (2000),  no. 1-2, 61--94.

\bibitem{Mo}{\sc L. J. Mordell}, The diophantine equation $y^2 -k = x^3$,  {\it Proc. 
Lond. Math. Soc. (2)}, {\bf 13}, (1913), 60-80. 

\bibitem{MoDio}{\sc L. J. Mordell}, {\bf Diophantine equations}, Pure and Applied Mathematics, Vol. 30 Academic Press, London-New York, 1969.

\bibitem{Na}{\sc J. Nakagawa}, On the relations among the class numbers of binary 
cubic forms, {\it Invent. Math.}  {\bf 134}  (1998),  no. 1, 101--138. 

\bibitem{Sh}{\sc T. Shintani}, On Dirichlet series whose coefficients are class 
numbers of integral binary cubic forms,  {\it J. Math. Soc. Japan}  {\bf 24}, (1972), 
132-188.

\bibitem{Sl-St}{\sc M. J. Slupinski and R. J. Stanton}, Symplectic geometry of 
Heisenberg graded Lie algebras, nearing completion.


\bibitem{Wr}{\sc D. J. Wright}, The adelic zeta function associated with the space of 
binary cubic forms, I: global theory,  {\it Math. Ann.} {\bf  270} (1985), 503--534.

\end{thebibliography}
\end{document}